\documentclass[10pt,twocolumn]{article} 
\usepackage{simpleConference}
\usepackage{times}
\usepackage{graphicx}
\usepackage{amssymb}
\usepackage{url,hyperref}

\usepackage[super]{natbib}

\usepackage{stackrel} 

\usepackage{amssymb,amsmath}

\usepackage[default]{lato}

\usepackage[utf8]{inputenc}
\usepackage[T1]{fontenc}
\usepackage{alphabeta}
\usepackage{stmaryrd} 
\usepackage{newunicodechar}
\newunicodechar{∀}{\ensuremath{\mathnormal{\forall}}}
\newunicodechar{→}{\ensuremath{\mathnormal{\to}}}
\newunicodechar{₀}{\ensuremath{{}_0}}
\newunicodechar{₁}{\ensuremath{{}_1}}
\newunicodechar{₂}{\ensuremath{{}_2}}
\newunicodechar{₃}{\ensuremath{{}_3}}
\newunicodechar{₄}{\ensuremath{{}_4}}
\newunicodechar{₅}{\ensuremath{{}_5}}
\newunicodechar{₆}{\ensuremath{{}_6}}
\newunicodechar{₇}{\ensuremath{{}_7}}
\newunicodechar{₈}{\ensuremath{{}_8}}
\newunicodechar{₉}{\ensuremath{{}_9}}
\newunicodechar{⁰}{\ensuremath{{}^0}}
\newunicodechar{¹}{\ensuremath{{}^1}}
\newunicodechar{²}{\ensuremath{{}^2}}
\newunicodechar{³}{\ensuremath{{}^3}}
\newunicodechar{⁴}{\ensuremath{{}^4}}
\newunicodechar{⁵}{\ensuremath{{}^5}}
\newunicodechar{⁶}{\ensuremath{{}^6}}
\newunicodechar{⁷}{\ensuremath{{}^7}}
\newunicodechar{⁸}{\ensuremath{{}^8}}
\newunicodechar{⁹}{\ensuremath{{}^9}}
\newunicodechar{⁺}{\ensuremath{{}^{+}}}
\newunicodechar{≟}{\ensuremath{\stackrel{?}{=}}}
\newunicodechar{∘}{\ensuremath{\circ}}
\newunicodechar{ᶠ}{\ensuremath{{}^f}}
\newunicodechar{ⁿ}{\ensuremath{{}^n}}
\newunicodechar{ᵐ}{\ensuremath{{}^m}}
\newunicodechar{ᵏ}{\ensuremath{{}^k}}
\newunicodechar{ⁱ}{\ensuremath{{}^i}}
\newunicodechar{ˢ}{\ensuremath{{}^s}}
\newunicodechar{ᵗ}{\ensuremath{{}^t}}
\newunicodechar{ᵃ}{\ensuremath{{}^a}}
\newunicodechar{ᵇ}{\ensuremath{{}^b}}
\newunicodechar{ᵛ}{\ensuremath{{}^v}}
\newunicodechar{ᵞ}{\ensuremath{{}^{\gamma}}}
\newunicodechar{⁻}{\ensuremath{{}^{-}}}
\newunicodechar{≡}{\ensuremath{\equiv}}
\newunicodechar{⊥}{\ensuremath{\perp}}
\newunicodechar{ℓ}{\ensuremath{\ell}}
\newunicodechar{∵}{\ensuremath{\because}}
\newunicodechar{∴}{\ensuremath{\therefore}}
\newunicodechar{∷}{\ensuremath{::}}
\newunicodechar{⸼}{\ensuremath{.}}
\newunicodechar{∈}{\ensuremath{\in}}
\newunicodechar{ℕ}{\ensuremath{\mathbb{N}}}
\newunicodechar{⦅}{\ensuremath{\llparenthesis}}
\newunicodechar{⦆}{\ensuremath{\rrparenthesis}}
\newunicodechar{•}{\ensuremath{\bullet}}

\IfFileExists{upquote.sty}{\usepackage{upquote}}{}
\makeatletter
\@ifundefined{KOMAClassName}{
  \IfFileExists{parskip.sty}{%
    \usepackage{parskip}
  }{
    \setlength{\parindent}{0pt}
    \setlength{\parskip}{6pt plus 2pt minus 1pt}}
}{
  \KOMAoptions{parskip=half}}
\makeatother
\usepackage{xcolor}
\usepackage{color}
\usepackage{fancyvrb}

\DefineVerbatimEnvironment{Highlighting}{Verbatim}{commandchars=\\\{\}}
\newenvironment{Shaded}{}{}

\newcommand{\CommentTok}[1]{\textcolor[rgb]{0.38,0.63,0.69}{\textit{#1}}}

\newcommand{\DataTypeTok}[1]{\textcolor[rgb]{0.56,0.13,0.00}{#1}}

\newcommand{\KeywordTok}[1]{\textcolor[rgb]{0.00,0.44,0.13}{\textbf{#1}}}
\newcommand{\NormalTok}[1]{#1}

\newcommand{\OtherTok}[1]{\textcolor[rgb]{0.00,0.44,0.13}{#1}}

\usepackage{bussproofs}
\usepackage{xargs} 
\usepackage{ifthen} 
\usepackage{amssymb}
\usepackage{stmaryrd} 
\usepackage{amsmath}
\newcommand{\myife}[4]{\ifthenelse{\equal{#1}{#2}}{#3}{#4}}

\newcommandx{\inference}[2][1=1, usedefault]{
	\myife{#1}{1}
	{\UnaryInfC{#2}}
	{\myife{#1}{2}
		{\BinaryInfC{#2}}
		{\myife{#1}{3}
			{\TrinaryInfC{#2}}
			{\myife{#1}{4}
				{\QuaternaryInfC{#2}}
				{\myife{#1}{5}
					{\QuinaryInfC{#2}}
					{}}}}}}

\newcommand{\Fin}[1]{\ensuremath{\mathbb{N}_{#1}}}
\newcommand{\Index}[1]{\ensuremath{\text{index } #1}}
\newcommand{\funII}[0]{\ensuremath{\text{ fun11 } }}
\newcommand{\funMI}[1]{\ensuremath{\text{funM1 } #1}}
\newcommand{\funIN}[1]{\ensuremath{\text{fun1N } #1}}
\newcommand{\funMN}[2]{\ensuremath{\text{funMN } #1 #2}}
\newcommand{\tuple}[1]{\ensuremath{\text{tuple } #1}}
\newcommand{\scalar}[0]{\ensuremath{\text{ scalar }}}
	
\usepackage{tikz-cd}
\usepackage{amsmath, amssymb, amsfonts}
\usepackage{mathtools}
\usepackage{tikz}
\newlength\figH
\newlength\figW
\usepackage{pgfplots}
\usepackage{subcaption}
\usepackage{float}
\usepackage{siunitx}
\usepackage{qtree}

\usepackage{enumitem}
\setitemize{noitemsep,topsep=10pt,parsep=0pt,partopsep=0pt}

\definecolor{green}{HTML}{FF7F00}
\definecolor{green}{rgb}{0.55, 0.71, 0.0}
\definecolor{green}{rgb}{0.0, 0.5, 0.0}

\newcommand*\diff{\mathop{}\!\mathrm{d}}

\begin{document}

\title{Encoding Electromagnetic Transformation Laws for Dimensional Reduction}

\author{Marcus Christian Lehmann, Mirsad Had{\v{z}}iefendi{\'{c}} \\
Albert Piwonski, Rolf Schuhmann\\
\\
Theoretische Elektrotechnik, Technische Universit{\"{a}}t Berlin\\
Berlin, Germany\\
\\
\today
}

\maketitle
\thispagestyle{empty}

\begin{abstract}
Electromagnetic phenomena are mathematically described by solutions of boundary value problems.
For exploiting symmetries of these boundary value problems in a way that is offered by techniques of dimensional reduction, it needs to be justified that the derivative in symmetry direction is constant or even vanishing. 
A generalized notion of symmetry can be defined with different directions at every point in space, as long as it is possible to exhibit unidirectional symmetry in some coordinate representation.
This can be achieved, e.g., when the symmetry direction is given by the direct construction out of a unidirectional symmetry via a coordinate transformation which poses a demand on the boundary value problem.
Coordinate independent formulations of boundary value problems do exist but turning that theory into practice demands a pedantic process of backtranslation to the computational notions. This becomes even more challenging when multiple chained transformations are necessary for propagating a symmetry. We try to fill this gap and present the more general, isolated problems of that translation.

Within this contribution, the partial derivative and the corresponding chain rule for multivariate calculus are investigated with respect to their encodability in computational terms. We target the layer above univariate calculus, but below tensor calculus.
\end{abstract}

\begin{quote}
This is the peer-reviewed but pre-published version of the following submitted article:

\textit{Lehmann MC, Hadžiefendić M, Piwonski A, Schuhmann R,
	Encoding Electromagnetic Transformation Laws for Dimensional Reduction.
	Int J Numer Model El. 2020;e2747.
	https://doi.org/10.1002/jnm.2747
}

which might be published in final form. This article may be used for non-commercial purposes in accordance with Wiley Terms and Conditions for Use of Self-Archived Versions.
\end{quote}




\newcommand\numberthis{\addtocounter{equation}{1}\tag{\theequation}}

\section{Introduction}\label{sec1}

There is a variety of different formulas for the transformation of vector components of fields and fluxes in classical electromagnetism.
When changing the coordinate-system the vector components need to be transformed because vector components quantify directions induced by the coordinate-system.
This results in a different matrix-transformation scheme, depending on the physical meaning of the vectors in question.
Different transformation properties of the objects considered in electromagnetic theories have been known for a long time\cite{Bossavit2001,VanDantzig1934}. They can be systematically formulated within tensor calculus at the cost of using antisymmetric tensors.
Representing electromagnetic objects with antisymmetric tensors leads to a high amount of combinatorics in tensor calculus, especially when resolving permutations.
The theory of differential forms provides a formalism to abstract over that. \par

Within the domain of computational electromagnetism exist several formalisms to represent physical entities with mathematical objects.
Most notably are quaternions, vectors, tensors, Clifford numbers and differential forms\citep{Bossavit2012}.
The borders of the theories for these objects, i.e. the precise number of logical laws belonging to each theory, are differently blurry.

\begin{figure}[h]
\begin{align*}
d_{\tau}E&=-\mathcal{L}_TB &\quad&&d_{\gamma}F^*E&=-\mathcal{L}_{\Gamma}F^*B\\
d_{\tau}H&=J+\mathcal{L}_TD&\quad&&d_{\gamma}F^*H&=F^*J+\mathcal{L}_{\Gamma}F^*D\\
d_{\tau}D&=\rho            &\quad\stackrel{\Large{\Leftrightarrow}}{\phantom{=}}\quad&&d_{\gamma}F^*D&=F^*\rho\\
d_{\tau}B&=0               &\quad&&d_{\gamma}F^*B&=0
\end{align*}
\caption{Formulation of a generic electromagnetic boundary value problem\cite{Raumonen2009} on the manifold $M$ with differential forms (left) and an equivalent boundary value problem on the manifold $N$ (right). Material laws and boundary conditions are omitted but they also follow the same pattern of transformation.}
\label{fig:raumonen}
\end{figure}

\hyphenation{theory}

When representing electromagnetic potentials, fields, fluxes and densities with \textit{differential forms}, the physical space is modeled by the notion of a \textit{manifold}.
An electromagnetic boundary value problem can be posed with the help of an \textit{observer structure}\cite{Raumonen2009} within the theory of differential forms.
For a differentiable manifold\footnote{see Sec. \ref{sec:electromagneticalcontext}} $M$, a smooth nonzero vector field $T$ on $M$ and a smooth one form $\tau$ on $M$ such that $\tau(T)=1$, the pair $(T,\tau)$ is called an \textit{observer structure}.
Using this observer structure $(T,\tau)$, the differential operators $d_\tau$ and $\mathcal{L}_T$ can be established.

A generic boundary value problem over one domain can be transformed to another domain by transforming the involved differential forms.
Boundary value problems are regarded \textit{equivalent}\cite{Raumonen2009} if they can be transformed into each other in that way.
If a boundary value problem is suitable for using techniques of dimensional reduction, then these techniques can be applied to all equivalent boundary value problems.
For an observer structure $(\tau,T)$ on a manifold $M$, an observer structure $(\gamma,\Gamma)$ on a manifold $N$, on these entities a transformation $F : N \rightarrow M$ and on the differential forms the induced transformation $F^*$, the two formulations of Fig. \ref{fig:raumonen} pose the same boundary value problem.
There is a very regular pattern present in these equations stating what needs to be done to transform a boundary value problem:
the differential forms of the original boundary value problem have to be transformed with the \textit{pullback} $F^*$ to appear in the transformed boundary value problem. An introduction into the calculus on manifolds and an introduction of differential forms can be found in the literature\cite{spivak71}.

On a machine, computations for solving a boundary value problem operate on \textit{number data} - the numbers that are stored within the machine's memory.
In the current formulation it might not be that obvious anymore how to convert the original number data into the number data for the transformed boundary value problem in terms of actual computations.
A confident implementation of a program benefits from an obvious description of the computation.
Therefore multiple formulations complement each other: for the computation, low-level matrix-operations and index-operations can directly be executed by the machine, but for deriving the computation, only the high level differential forms statements can be overviewed. We are convinced, that high level abstractions as in Fig. \ref{fig:raumonen} pay off in the most beneficial way only, when stacked on top of a layer providing
\begin{itemize}
	\item[a)] a good abstraction to provide coordinate transformation rules in terms of matrix-based or just general computation schemes for a given tensorial formulation and
	\item[b)] a good abstraction to incorporate combinatorial notions, especially the enumeration of permutations, which enables the reasoning on a level of differential forms to be automatically transferred into a tensorial representation.
\end{itemize}
The purpose of an implementation is to put the machine into a state that is most efficient for processing all necessary computations of a numerical scheme which solves the boundary value problem.
Abstractions help in organizing the implementation but should not prevent to use the machine in it's most efficient way.
Therefore most abstractions are usually stripped before a cost-intensive computational task is started.
They should only allow to produce an efficient computational scheme on spot in some form that is available on the machine: matrix or parallel or other kinds of efficient computational primitives.
It is important to emphasize that the corresponding raw number data does not need to change for every transformation process.

For the first part a), i.e. the generation of transformation rules, in this paper, we show a way to realize such a layer which is independent of the actual function representation. The second part b) is motivated in Sec. \ref{sec6} and not treated in this paper.

The rest of this paper is organized as follows:
Since this is an interdisciplinary topic, we will give some necessary context for readers from different domains.
This context is tailored to an implementation on a machine.
In Sec. \ref{sec:electromagneticalcontext} we introduce the notion of the \textit{frame bundle} and \textit{associated bundles} to the frame bundle from \textit{bundle theory} in order to define \textit{geometric quantities} and give the general transformation rule for $(p,q)$-tensor $\omega$-densities.
In Sec. \ref{sec:softwarecontext} we introduce shape functions and degrees of freedom of the finite element method in the context of differential forms.
An introduction to the untyped $\lambda$-calculus is given in Sec. \ref{sec:computationalcontext}.
A definition of the partial derivative in terms of a univariate derivative is given in Sec. \ref{sec:problem}.
In Sec. \ref{sec:theory} the untyped lambda calculus is augmented with axioms for a typed variant for the purpose of expressing the calculations of the previous section.
In Sec. \ref{sec:software} we give a guideline on how this augmented lambda calculus can be applied in a software project.


\subsection{Electromagnetical Context}\label{sec:electromagneticalcontext}

\begin{figure}[h]
	\centering
	\def\svgwidth{0.9\linewidth}
	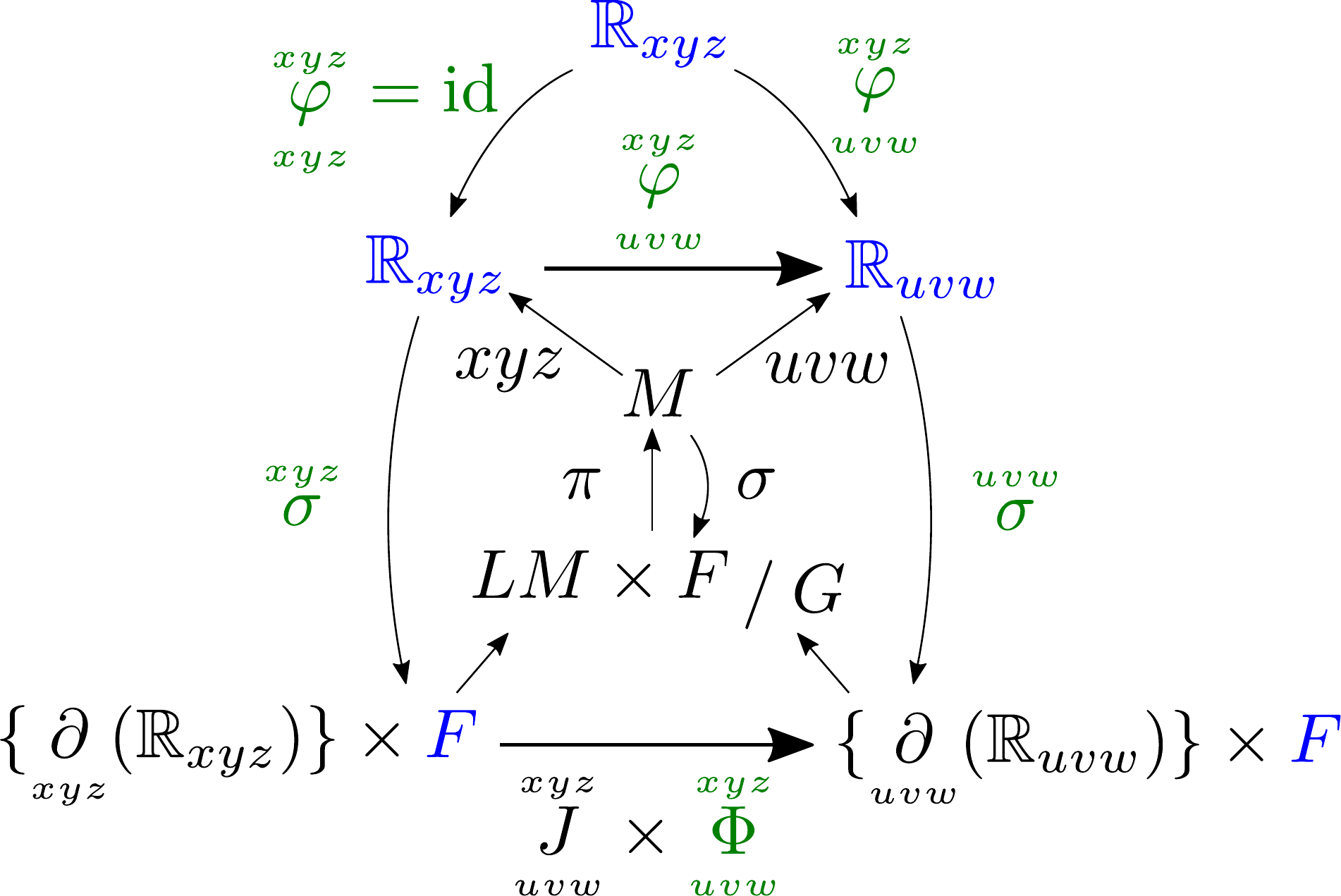
	\caption{Objects involved in a \textit{covariant} treatment of associated fibre bundle. All arrows in this diagram denote functions. The name of the function is written next to its arrow. At the beginning of an arrow is the domain space and at the end of an arrow is the co-domain space of the corresponding function. A product of spaces is denoted by $\times$ and similarly the parallel composition of two functions operating on these product spaces is also denoted with $\times$.}
	\label{fig:diagram}
\end{figure}

%

Electromagnetic theory is concerned about the spatial and temporal relation of different physical quantities such as potentials, forces, fluxes and densities\cite{Tonti2013}.
These are often grasped with respect to a coordinate system and its \textit{coordinate-induced directions}. A base for all directions at a point is called a \textit{frame}\cite{hehl03}.
The coordinate system is called a \textit{chart} and it is modeled as a continuous mapping from \textit{points} $p$ of a topological space to \textit{their} coordinates within $\mathbb{R}^n$. A collection of such systems is called an \textit{atlas} and if the whole space can be covered by overlapping charts into $\mathbb{R}^n$ for the same $n$ it is called \textit{locally euclidean}. A \textit{topological manifold} is defined by additionally demanding the \textit{Hausdorff}\cite{Raumonen2009} property and the space being \textit{second countable}\cite{Raumonen2009}. If chart transition functions of an atlas are arbitrarily often differentiable, the atlas is called \textit{smooth}. Two atlases over $M$ are \textit{smoothly equivalent} when their union is also a smooth atlas over $M$. An equivalence class of smoothly equivalent atlases over $M$ is called a \textit{smooth structure}. 
A topological manifold $M$ endowed with a smooth structure is called a \textit{smooth manifold}.
Analogously, regarding only $k$-times differentiable chart transition functions for $k>0$ leads to a \textit{differentiable structure}.
A topological manifold endowed with a differentiable structure is called a \textit{differentiable manifold}\cite{Raumonen2009}.
On the differentiable manifold, we defined the observer structure that was necessary to establish the differential operators for the boundary value problem in the introduction.
A far-reaching introduction on this topic can be found in the literature\cite{baez94}\cite{hehl03}.
Having two manifolds, one called \textit{total space} $E$ and one called base space $M$, and a continuous surjective function $\pi:E\to M$, the tuple $(E, \pi, M)$ is called a \textit{bundle} of manifolds.
Here the preimage $\text{preim}_\pi(p)$ of a point $p\in M$ with respect to $\pi$ is called \textit{fibre at the point $p$}, denoted by $F_p$. If fibres of all points are homeomorphic to some manifold $F$, the bundle is called a \textit{fibre bundle with typical fibre $F$}.
Globally over the manifold, tensor \textit{fields}, vector \textit{fields} and differential forms are considered. They are modeled as sections of some fibre bundle where locally at a point tensors, vectors and co-vectors of some algebra are considered.
One specifically important algebra for that purpose is the \textit{exterior algebra} of local co-vectors from global differential forms.

When answering \textit{``why exterior differential forms''} are useful as a formalism for the modeling of electromagnetic laws, some authors\cite{hehl03} justify this with \textit{``the alternating algebraic structure of integrands that gave rise to the development of exterior algebra and calculus which is becoming more and more recognized as a powerful tool in mathematical physics''}\footnote{p. 21, Hehl\cite{hehl03}}.
Further, we will make use of the generalization of a tensor, the \textit{geometric quantity}, being \textit{``defined by the action of the general linear group on a certain set of elements''}\footnote{p. 27, ibd.}. Examples are tensor-valued differential forms and twisted tensors.

Electromagnetism as a physical effect does not depend on a chosen coordinate system, which - as a property - is called \textit{general covariance}. In our new wording, a geometric quantity at some point should not depend on a chosen frame.
Now, a standard mathematical conjuring trick in order to avoid an arbitrary choice, is to attach that choice to the objects in question - to attach the chosen frame to the quantity in our case.

The theory of associated fibre bundles can describe different kinds of fibre bundles that fulfill an equivalence relation.
This equivalence relation is expressed by the means of a Lie group $G$ with respect to some $G$-bundle. A $G$-bundle will be introduced straightaway.

For a manifold $M$, a cover $\{U_\alpha\}$ of $M$ by open sets, a vector space $V$ and a group representation $\rho:G\to GL(V)$ of $G$ on $V$, where $GL(V)$ is the general linear group over $V$, it is possible to obtain\cite{baez94} a vector bundle $(E,\pi,M)$ using transition functions $g_{\alpha\beta}:U_\alpha\cap U_\beta\to G$.
This is only possible when the \textit{compatibility conditions}\cite{baez94}
\begin{align*}
g_{\alpha\alpha}&=1&&\quad\text{on }U_\alpha\\
g_{\alpha\beta}g_{\beta\gamma}g_{\gamma\alpha}&=1&&\quad\text{on }U_\alpha\cap U_\beta\cap U_\gamma.
\end{align*}
are fulfilled. It is done\footnote{adapted from p. 212-214, Baez\cite{baez94}} by partitioning the disjoint union $\cup_\alpha\left(U_\alpha\times V\right)$ with the equivalence relation
$$(p,v)\sim (p,v')\;\Leftrightarrow\;v=\rho(g_{\alpha\beta}(p))\,v'$$
to obtain the base space $E=\cup_\alpha\left(U_\alpha\times V\right)/{\sim}$.
A vector bundle obtained in this way is called a $G$-bundle and $F$ is the \textit{standard fibre}\cite{baez94}.
When $V=G$, the $G$-bundle is called \textit{principal}\cite{auchmann14}.
The frame bundle $LM$ is a principal $G$-bundle where $G$ is the general linear group.

The fibre bundle of frames over a smooth manifold $M$ is called \textit{frame bundle} and denoted by $LM$. At a point $p$ for some chosen frame $e\in L_pM$ and some geometric quantity $f\in F_p$ from the fibre $F_p$, we regard the tuple $(e,f)$ as one representation of a geometric quantity at that point.
But one could choose another frame $e'$ which is done by choosing another chart, or coordinate system, and then the new frame is related $\sim$ to the old one by means of the Jacobian $J$ at that point $p$:
$$(e,f)_p\sim(e\triangleleft J,J^{-1}\triangleright f).$$
Here $\triangleleft$ is a right action on the frames
$$(\vec{e}_1,...,\vec{e}_d)\triangleleft J:=(J^m_1\vec{e}_m,...,J^m_d\vec{e}_m)$$
and $\triangleright$ is a left action on the fibre
\begin{align*}
&{(J^{-1}\triangleright f)^{i_1...i_p}}_{j_1...j_q}:=\\
&\quad(\det J)^{\omega}\,{f^{k_1...k_p}}_{l_1...l_q}{(J^{-1})^{i_1}}_{k_1}...{(J^{-1})^{i_p}}_{k_q}{J^{l_1}}_{j_1}...{J^{l_p}}_{j_q}
\numberthis\label{eq:generaltransformationlaw}
\end{align*}
for all, so-called $(p,q)$-tensor $\omega$-densities. Both definitions make use of a \textit{sum convention}, summing over all equal indices.

The Jacobian $J$ is an element of the general linear group from the definition of a geometric quantity.
That makes $(p,q)$-tensor $\omega$-densities a special case of geometric quantities.
Not one tuple $(e,f)$ with one chosen frame, but the equivalence class of all tuples that can be related $\sim$ to each other with some $J$ makes a \textit{value} of a geometric quantity at a point $p$.
This is expressed in the inner part of Fig. \ref{fig:diagram}: by taking the space $LM\times F$ but partitioning it $(LM\times F)/G$ with respect to a group $G$. This group is the general linear group in our case.
The new fibre bundle $(LM\times F)/G$ is called to be \textit{associated} to the \textit{principal} $G$-bundle $LM$.

For the implementation on a machine, we are most likely not to work with a point $p\in M$ but with its coordinates $xyz(p)\in \mathbb{R}^3$.
These coordinates are with respect to some chart $xyz:M\to\mathbb{R}^3$ where $\mathbb{R}^3$ is denoted as $\mathbb{R}_{xyz}$ in Fig. \ref{fig:diagram}.
That chart $xyz$ induces directions and in particular one concrete frame $\stackbin[xyz]{}{\partial}(\mathbb{R}_{xyz})$ at each point $p$. To this frame, there corresponds exactly one number data from the fibre $F$ so as programmers we consider the \textit{chart representation} $\stackrel{xyz}{\sigma}$ of a section $\sigma$ in the implementation. When changing charts, the chart transition map $\stackbin[uvw]{xyz}{\varphi}$ transforming the coordinates has an analogue $\stackbin[uvw]{xyz}{J}\times\stackbin[uvw]{xyz}{\Phi}$ operating on the field representations. Generally we have
$$\Phi(f):=J^{-1}\triangleright f.$$
You can find an interpreted version of Fig. \ref{fig:diagram} in Sec. \ref{sec:software} as Fig. \ref{fig:interpreteddiagram}.

Usually one only represents the blue bits of Fig. \ref{fig:diagram}as data and the green bits of Fig. \ref{fig:diagram}as computations in an implementation. The remaining black bits might be treated in an \textit{opaque} way. This is a technique in creating programming interfaces where objects are exposed via references which are of a defined reference type. That reference is called \textit{opaque} when referring to unexposed or even undefined data while the representation of the reference itself is known\footnote{References are usually represented as integers indicating an item over some universe, most commonly an enumeration of cases or an address of the machine's memory.}. Even though the black bits are not themselves represented in an implementation as number data, their rules of operation are a candidate for entering an implementation as rules of opaque references.
Opaque references can be used to restrict the usage of operations on number data to valid cases. The amount and flexibility of expressible restrictions for that purpose is a property of the targeted programming language.
In Sec. \ref{sec:theory} we make use of a \textit{dependently typed}\cite{norell2009} programming language that offers high flexibility in expressing restrictions to increase confidence in our approach.
In our performance critical code we make use of a \textit{deterministic just-in-time-compiled} programming language that offers partial recompilation and high flexibility for code-generation to help putting the machine into its most efficient state for a computation.

The theory of associated fibre bundles of the frame bundle $LM$ provides a notion of $(p,q)$-tensor $\omega$-density. This notion is sufficient to express all electromagnetic quantities of interest and they share a single transformation law \ref{eq:generaltransformationlaw}.
The transformation $F^*$ of Fig. \ref{fig:raumonen} used to transform an electromagnetic boundary value problem follows the rules of $\Phi$ from Fig. \ref{fig:diagram}: having a single explicit definition for all the various $(p,q,\omega)$ transformations, makes this theory very promising as a starting point for an implementation in a software project.
Furthermore, it is observable that such software will heavily rely on the correct evaluation of Jacobians at the right coordinates of possibly chained transformations.
That is the reason why we are so interested in a very solid foundation of encoding partial derivatives and their chain rule in Sec. \ref{sec:problem} and \ref{sec:theory}.

\subsection{Software Context}\label{sec:softwarecontext}

There is much within computational electromagnetism that counts as \textit{software}. 
This community has a history of incorporating guidance from \textit{mathematical structure} of the electromagnetic theory into the development of \textit{consistent} and \textit{stable} numerical methods.

%
%
%
%
%
%

When speaking of numerical software, this paper focuses on the finite element method which is a \textit{Galerkin method} that can be expressed in terms of the \textit{Finite Element Exterior Calculus}\cite{arnold10}.
There are a lot of common mistakes leading to wrong solutions of the finite element method. The reason of failure is often not that obvious\cite{arnold10}.
This non-obviousness is mirrored in an extensive development of Galerkin methods, and in particular the finite element method, within past decades.




For a \textit{numerical} consideration, i.e., for the purpose of establishing proven guarantees of certain \textit{errors}, a numerical method is abstractly modeled using an abstract Hilbert space $V$.
It is assumed that the numerical problem can be expressed using a bounded bilinear form $B:V\times V\to\mathbb{R}$ and a bounded linear form $F:V\to\mathbb{R}$ as
\begin{equation}
\text{find }u\in V\text{ such that }\forall\,v\in V\,.\,B(u,v)=F(v).
\end{equation}
The problem is called \textit{well-posed} if an unique solution $u$ exists and the \textit{solution mapping} $F\mapsto u$ is bounded again.
Using that formulation, a Galerkin method is characterized by a family of \textit{finite dimensional, normed} spaces $V_h$ indexed by parameter $h$ that in \textit{some sense}\footnote{The spaces $V_h$ do not necessarily have to be subspaces of $V$.} approximate $V$.
The Galerkin method for that family of spaces $V_h$, a bilinear form $B_h:V_h\times V_h\to\mathbb{R}$ and a linear form $F_h:V_h\to\mathbb{R}$ poses
\begin{equation}
\text{find }u_h\in V_h\text{ such that }\forall\,v\in V_h\,.\,B_h(u_h,v)=F_h(v).
\end{equation}
It is desired to prove that the property of \guillemotleft $V_h$ approximating in \textit{some sense} $V$ as $h$ \textit{advances}\guillemotright{} is conveyed to \guillemotleft $u_h$ approximating in \textit{some sense} $u$ as $h$ advances\guillemotright.
The \textit{finite} element method is a Galerkin method where elements of the basis of $V_h$ have \textit{finite} support, i.e. they are nonzero only on a small part of the considered domain.


A construction of bases for a family $V_h$ of spaces can be proven to be consistent and stable when used in a Galerkin method.
The Finite Element Exterior Calculus provides constructions of classes of finite element bases whose Galerkin methods were proven to be consistent and stable\cite{arnold09}\cite{arnold10}.
This was done utilizing notions from differential geometry and algebraic topology in order to develop methods for error analysis.
It is necessary to do this within a functional analytic setting, because a notion of \textit{approximating} and \textit{error} and therefore \textit{consistency} and \textit{stability} of a numerical method, do ultimately origin here.





We consider the explicit construction\cite{arnold09} of two families of explicit local bases.
Here the approach was \textit{``not trying to find hierarchical bases, but rather [...] generalize the explicit Bernstein basis''}\cite{arnold09}.
Where it is easy to give a \textit{spanning set} of polynomials with meeting requirements, it is much harder\cite{arnold09} to provide a basis of \textit{linearly independent} polynomials.

A key insight is to decompose this construction of base elements and define the polynomial base in terms of smaller \textit{shape functions}.
Multiple adjacent of those shape functions are \textit{recombined} into one base element by enforcing \textit{proper interelement continuity} conditions\cite{arnold09}.
This approach is sometimes called an \textit{assembly}\cite{arnold09}~\footnote{The buildup of a matrix for the resulting discrete linear system is also called assembly process.}.
With its interelement continuity conditions, this process is the reason why multiple shape functions from adjacent pieces of a domain share the same \textit{degree of freedom}.
The presented\cite{arnold09} assembly process for the construction of basis functions, is \textit{``a straightforward consequence of the geometric decomposition of the finite element spaces''}\cite{arnold09}.



For both families of shape function spaces for each \textit{simplex} $T$ and each subsimplex, sometimes called \textit{face}, $f$ with $r\geq 1$ and $0\leq k\leq d$ and $d=\dim f \geq k$ there is a shape function space and there are \textit{degrees of freedom}.
One is the shape function space of polynomial differential forms $\mathcal{P}_{r}\Lambda^k(T)$ with corresponding degrees of freedom
$$u\mapsto \int_f(\text{tr}_f u)\wedge q \quad:\quad \mathcal{P}_{r}\Lambda^k(T)\to\mathbb{R},$$
where $q$ is from $\mathcal{P}^{-}_{r+k-d}\Lambda^{d-k}(f)$ and $\wedge$ is the exterior product. For differential forms, the trace operation $\text{tr}$ is the pullback $\iota_{f,T}^*$ of the inclusion $\iota_{f,T}:f\to T$
$$\text{tr}_f = \iota_{f,T}^*.$$
The other one is the shape function space of polynomial differential forms $\mathcal{P}^{-}_{r}\Lambda^k(T)$ with corresponding degrees of freedom given by
$$u\mapsto \int_f(\text{tr}_f u)\wedge q \quad:\quad \mathcal{P}^{-}_{r}\Lambda^k(T)\to\mathbb{R},$$
where $q$ is from  $\mathcal{P}_{r+k-d-1}\Lambda^{d-k}(f)$.

The construction of the base elements for a shape function space is \textit{``somehow a complicated business''}\footnote{\textit{Finite element exterior calculus} NSF/CBMS course, ICERM, 2012; lecture~9, min.~40 sec.~10.} and provided\cite{arnold09} in terms of:
\begin{itemize}
	\item a simplicial complex
	\item taking the set of subsimplices of a given simplex
	\item restriction maps from a simplex to one of it's subsimplices and inclusion maps the other way around
	\item barycentric coordinates
	\item the exterior derivative of barycentric coordinates
	\item piecewise polynomial differential forms
	\item the pullback of polynomial differential forms along a restriction map
	\item multi-indices
	\item the index-set associated to a face of the simplicial complex
	\item the set of all order preserving maps of indices
	\item taking the \textit{support} of a multi-index
	\item taking the \textit{range} of an order preserving map.
\end{itemize}
For every item on this list, we will probably have some correspondence within an implementation for a machine.
A simplicial complex is usually given by a \textit{mesh}.
It is mostly stored in two separate parts. One part is an abstract simplicial complex consisting just of the combinatorial information which is sometimes called the \textit{mesh topology}.
The other part is additional data which can be used to create homeomorphisms from the standard simplices to the given ones.
This data forms a \textit{parametrization} and provides barycentric coordinates.
In the case of a simplicial complex, this data might just contain vertex coordinates, but it becomes more interesting for \textit{curved} cells.

Multi-indices and order-preserving index maps are rooted in the combinatorial domain. Their representation in an implementation might be exploited in a clever way.
The most intriguing correspondence we think is the one of polynomial differential forms and their pullbacks.
These can be resolved within a pen and paper computation\footnote{A limited amount is shown in a table from the original paper\cite{arnold09} and the various families of bases are implemented within the FENICS\cite{fenics} project. Reproducing the bases from their paper required us some amount of bookkeeping. } and, then, the resulting polynomials can be implemented very carefully. But it also seems reasonable to formulate the whole construction of a shape function element within a programming language.
One of our goals within this paper is investigating how to do so in an appropriate way.

This approach essentially lifts the implementation to a meta-level. Previously, as programmers, we were seeking an implementation to perform a numerical computation in the most efficient way for a given machine.
Now, we have to program an implementation that is able to produce another, more concrete, implementation which in turn is able to perform the numerical computation in the most efficient way for a given machine.
The efficiency of this meta-implementation is usually not critical for the efficiency of the resulting implementation.

One obvious technique is to generate source code of an implementation with the meta-implementation.
The programming language of the meta-implementation does not need to be the same as the programming language for the targeted implementation.
Most programming languages offer meta-programming constructs to generate computations and data structures and the expression of constraints that are to be checked during this generation.
These constraints are used to restrict the argument's domain of a meta-computation.
The templating system\cite{cpp17} of the C++ programming language is a very popular choice in the community of computational electromagnetism\cite{pellikka2013}.
This might be partly because it allows to use the same language for the implementation and the meta-implementation.
While this choice of programming language helps the programmer in putting the machine into its most efficient state, it offers limited flexibility in expressing logical constraints for the valid application of meta-computations.
Therefore, the expression of algebraic rules from a construction of finite element bases might only be partially incorporated.
When seeking for confidence, it is critical to be able to express all rules that are needed to be confident of.
These rules are expressed in the programming-language of the meta-implementation in order to have them checked automatically.
We might even claim that the usefulness of checking rules critically depends on the completeness, or coverage, of the rules regarding \textit{all possible cases}.
Putting it in another way: we claim that
\begin{itemize}
    \item achieving high efficiency is the biggest challenge when programming the implementation, whereas
    \item achieving high validity is the biggest challenge when programming the meta-implementation.
\end{itemize}
That is why we advocate the use of a programming language with a checking mechanism for \textit{dependently typed}\cite{norell2009} expressions to formulate constaints.
For the meta-program this offers a chance to express all algebraic rules completely.
This is relevant because meta-implementation techniques seem to become more and more unavoidable in modern high performance computing.

\subsection{Computational Context}\label{sec:computationalcontext}

Within this contribution, the partial derivative and its corresponding chain rule for multivariate functions will be investigated with respect to their encodability in \textit{computational terms}.
A functional analytic setting is very powerful for an analysis of problems related to partial differential equations.
In this paper, we will treat the operation of taking the derivative of a univariate function in a more \textit{synthetic} way.
The derivative operation will be embedded into a more general context of computation where some basic properties become assumptions of that embedded derivative operation.

In this paper, the understanding of computational terms is backed by lambda-calculus ($\lambda$-calculus) which serves as a model, or definition, of \textit{effectively calculable} functions.
That calculus was originally developed by A. Church in 1936\cite{church36} and we will follow a modern treatise\cite{barendregt85} of the resulting findings.
We will take a \textit{type free}\footnote{In our reference\cite{barendregt85} this is in preparation for explaining Curry-style and Church-style $\lambda$-calculi} $\lambda$-calculus that is extended in Sec. \ref{sec:theory} to a typed variant.
The type-free $\lambda$-calculus is constituted by a set $\Lambda$ of $\lambda$-terms built up from an infinite set of \textit{variables} $V=\{v,v',v'',...\}$ using \textit{application} and \textit{function abstraction}:
\begin{align*}{}
x\in V&\quad\implies\quad x\in\Lambda,&{}\\
M,N\in\Lambda&\quad\implies\quad (M\,N)\in\Lambda,&{}\\
M\in\Lambda,x\in V&\quad\implies\quad (\lambda x\, M)\in\Lambda.&{}
\end{align*}
We choose the convention to suppress the outermost parenthesis in $(\lambda x\, M)$ when it is unambiguous and to add a separating dot inbetween $x$ and $M$ resulting in $\lambda x\,.\, M$.
On these terms, an operation of \textit{substituting} $N$ for the free occurences of $x$ in $M$ can be defined and is denoted by $M[x:=N]$.

Furthermore there are binary relations for $\eta$-reduction, $\alpha$-conversion and $\beta$-reduction reading from left to right:
\begin{align}{}
\lambda\, x \,.\, f\,x&\;\equiv_{\eta}\;f\label{eq:Lambda-eta}\\
\lambda\, x \,.\, f\, x&\;\equiv_{\alpha}\; \lambda\, y \,.\, f\,y\label{eq:Lambda-alpha}\\
\left(\lambda\, y \,.\, M\right)\, x &\;\equiv_{\beta}\;M[y:=x]\label{eq:Lambda-beta}.
\end{align}
We will regard two terms as computationally equivalent if they can be related to each other in these ways.
Therefore a symmetric \mbox{$\equiv$-symbol} is already present here, although $\eta$-reduction, $\alpha$-conversion and $\beta$-reduction are defined as operations from the left hand side to the right hand side in the previous listing. The usage of $\lambda$-calculus will be elaborated in more detail in Sec.~\ref{sec:problem} and put in a more rigorous setting in \mbox{Sec. \ref{sec:theory}}.

An introductory survey, reaching out to the techniques\footnote{Encoding of de Bruijn indices\cite{debruijn72} can be regarded with respect to that interpretation, but does not need to. A further treatment of defining equivalences and proving their preservation, makes the \textit{encoding} of propositions as types indispensable. Therefore propositions become common objects to handle in a formal implementation.} used in Sec. \ref{sec:theory}, can be found in the literature\cite{wadler15} as the \textit{propositions as types} paradigm.
This paradigm pictures the development from $\lambda$-calculus to the \textit{proofs-as-programs} and \textit{propositions-as-types} interpretation through one of the most prominent developments within theoretical computer science: the Curry-Howard correspondence.

The rules of the black bits of Fig. \ref{fig:diagram}, mentioned in Sec. \ref{sec:electromagneticalcontext}, lead to constraints for restricting an argument's domain of a meta-computation, mentioned in Sec. \ref{sec:softwarecontext}.
These rules can be formulated as propositions of objects within the theory of differential forms.
Analogously, the constraints for restricting an argument's domain of a meta-computation can be formulated as propositions of an object in the programming language.
A compiled programming language is able to perform compile-time checks based on type-equations formulated in the programming language's type system.
Expressing a proposition as type-equations, and having these type-equations checked, or affirmed, conveys the affirmation of the type-equations to an affirmation of the proposition.
Therefore, using the type system of a programming language for the purpose of establishing validity of a meta-implementation depends on the translatability of propositions into equivalent type equations - equivalent with respect to affirmation. This orientation towards translatability and high validity differs from the algorithmically-oriented approach of computer algebra systems.

To support the theory in Sec. \ref{sec:theory} we have used a programming language that is developed precisely for the purpose of establishing a translatability of propositions into equivalent type equations.
This choice seems to offer the best chance of being able to express all algebraic rules of a construction of finite element shape functions completely.
But that is an outlook. In this paper we propose an encoding of the partial derivative in $\lambda$-calculus as a foundation of a system of algebraic rules.
That foundation is tailored towards an application in numerical methods, especially the finite element method.

\section{Problem and Application}\label{sec:problem}

\subsection{Varying Syntax and Semantics}\label{sec2b}

We introduced lambda calculus in the standard notation, which is also the notation we use in our implementation later-on.
But for this \textit{modeling} part, we switch to a barred-arrow $\mapsto$ notation since it resembles the standard notation of the electromagnetic theory.

References we give here for notation, might be a bit \textit{picked}, and it is, of course, a matter of taste. But it is this notation that illustrates day-to-day problems when working within the electromagnetic theory.

In one reference\cite{bossavit98} from the domain of computational electromagnetism, A. Bossavit argues\footnote{A.1.9 A notation for functions} about a notation for functions.
The given argument is to advertise using an arrow-symbol\footnote{A. Bossavit uses a straight-arrow $\to$ for both, function abstractions and function types, whereas we would argue to use the barred-arrow $\mapsto$ for function abstractions and $\to$ for function types.} in order to better emphasize a distinction of functions and \textit{expressions}. It is recommended to denote a function $f$ of the expression $x^2 + 2x + 1$ as
$$f = x\to x^2 + 2x + 1$$
or rather
\begin{equation}
f(x)\,:=\, x^2 + 2x + 1.
\label{eq:codomaindefinition}
\end{equation}
Where it is stressed that using just $=$ in the second case, could be interpreted as an equality instead of a definition. This is accompanied with an example of a differential operator, helping to resolve some \textit{``ambiguity as to which gradient, with respect to x or to y, we mean''}, making $x$ \textit{the parameter} and $y$ \textit{the variable}, both of which are vectors:
$$\text{grad}\left(y\to \frac{1}{|x-y|}\right)=y\to\frac{x-y}{|x-y|^3}.$$
Those differential operators act on function objects, and their notation might be borrowed from the notation of \textit{higher order functions} in programming.
A reference to programming, and especially to $\lambda$-calculus is already drawn in that reference.

In the same way that higher order functions, or functional programming in general, are known to have some steep learning curve to overcome, similar applies here.
This might be, why usually in engineering a \textit{codomain-focused style of notation}, as in (\ref{eq:codomaindefinition}), is preferred.

We think that some confusion arises by taking expressions and not functions as the dominant objects in calculus.
For instance, the Mathematica\cite{mathematica} programming language follows an expression focused approach.

Speaking about expressions, coincidentally also another example\cite{martinlof85} is given by P. Martin-L\"of, although he was not up for differential calculus and used it as a mere example for \textit{forms of expressions}:
\begin{quotation}
\noindent The expressions [...] are formed from variables
$$x,y,z$$
by means of various forms of expression
$$(F\,x_1,...,x_n)(a_1,...,a_m).$$
In an expression of such a form, not all of the variables $x_1,...,x_n$ need become bound in all of the parts $a_1,...,a_m$.
Thus, for each form of expression, it must be laid down what variables become bound in what parts. For example,
$$\int_a^b f\,\diff x$$
is a form of expression $(I\,x)(a,b,f)$ with $m=3$ and \mbox{$n=1$}, which binds all free occurrences of the single variable $x$ in the third part $f$, and
$$\frac{\diff f}{\diff x}(a)$$
is a form of expression $(D\,x)(a,f)$, with $m=2$ and \mbox{$n=1$}, which binds all free occurrences of the variable $x$ in the second part $f$.
\end{quotation}
\pagebreak
Furthermore, M. Spivak is\cite{spivak71} denoting the multivariate and univariate chain rule by
$$D(g\circ f)(a)=Dg(f(a))\circ Df(a)$$
$$(g\circ f)'(a)=g'(f(a))\cdot f'(a)$$
and introduces the $i$-th partial derivative of $f$ at $a$ as $D_if(a)$. He alludes, that the partial derivative is the ordinary derivative of a certain function, e.g. if \mbox{$g(x)=f(a^1,...,x,...,a^n)$} then 
\begin{equation}
D_if(a)=g'(a^i).\label{eq:ordinaryderivativeofcertainfunction}
\end{equation}
Notation is briefly discussed and it is mentioned that $D_1f(u,v,w)$ resolves the usage of a notation like
$$\left.\frac{\partial f(x,y,z)}{\partial x}\right|_{(x,y,z)=(u,v,w)}\;\text{ or }\;\frac{\partial f(x,y,z)}{\partial x}(u,v,w).$$
Another issue is framed by E. Tonti who also dedicates a chapter\cite{Tonti2013} to revise terminology.
There, \textit{many kinds of equality} are illuminated. He proposes five different such equalities to be suitable for the purpose of explaining electrodynamics instead of just using a single $=$ for all of them:
\begingroup
\begin{align*}
\stackrel{\text{def}}{=} & \text{\scriptsize~definition }   & H                                     & \stackrel{\text{def}}{=} U+p\,V, \;\text{\scriptsize~definition of enthalpy} \\
\equiv                   & \text{\scriptsize~identity }     & a^2-b^2                               & \equiv (a+b)(a-b),               \;\text{\scriptsize~for all }a\text{\scriptsize~and }b \\
=                        & \text{\scriptsize~equation }     & 3x^2 - 2x                             & = 5,                             \;\text{\scriptsize~the variable }x\text{\scriptsize~is unknown} \\
\stackrel{\text{mat}}{=} & \text{\scriptsize~material law } & V                                     & \stackrel{\text{mat}}{=} R\,I,   \;\text{\scriptsize~Ohm's law} \\
\stackrel{\text{law}}{=} & \text{\scriptsize~general law }  & \partial_t\rho+\text{div}\,\mathbf{J} & \stackrel{\text{law}}{=} 0,      \;\text{\scriptsize~conservation law}.
\end{align*}
\endgroup
This issue could also summarized by arguing that \textit{``the fragment of mathematical symbolese available to most calculus students has only one verb, `=' ''}\cite{thurston94}. \textit{``That’s why students use it when they’re in need of a verb.''}\cite{thurston94} In general, there is \textit{``a list of different ways of thinking about or conceiving of the derivative''}\cite{thurston94} of a function instead of a single way to do so. These all make their appearance at some point when studying electromagnetism.

We might summarize that these authors propose an expressive notation for what kind of statement is expressed by $=$, maybe even which kind of objects it relates, and how the variables of an expression are quantified.


Rather than giving a meaning of what \textit{the} univariate derivative is, we treat it synthetically and collect the few properties necessary for introducing a partial derivative on top. In order to resolve the various notations, we have chosen to resemble $\lambda$-calculus.


To support multiple interpretations, we chose to explain our usage of $\lambda$-calculus with a changed notation.
From our experience this better resembles day-to-day notation in computational electromagnetism but is close enough to follow notation of a formalization later-on in Sec.~\ref{sec:theory}.
This choice is made to support readers that are not immediately implementing such $\lambda$-calculus but still want to gain some insights about the partial derivative.


\subsection{Yet another notation for functions}

Our aim is to connect more high level theories, such as tensor calculus and differential forms to more low level theories, such as multivariate calculus and $\lambda$-calculus.
With tensors and differential forms it is possible in a tractable way to express sound notions of invariant properties and differentials.
In multivariate calculus and $\lambda$-calculus it is possible in a tractable way to express sound notions of an univariate derivative and computations.
After such a connection is made, representations and implementations that arguably behave in a way respecting these notions need to be given. Doing so should contribute to the discussion about how higher level \textit{representations of physical entities} can be encoded in a program.

We start with the assumption of a given univariate derivative operation~${\color{red}\mathbf{'}}$ that for a given univariate function representation $f$ can compute the univariate function representation of the derivative of that function $f{\color{red}\mathbf{'}}$. For the computational description, we make use of an untyped, simplified $\lambda$-calculus as introduced in Sec \ref{sec:computationalcontext}. Instead of $\lambda x. f x$, we denote function abstraction by $x\mapsto f(x)$ to better resemble day-to-day notation. We emphasize that only the following rules are used and it does not matter if you do not know $\lambda$-calculus yet, if you can familiarize yourself with these four \textit{computational equivalences} (\ref{eq:lambda-eta}-\ref{eq:lambda-comp}) that are already in use in engineering mathematics and denoted by~$\equiv$ here. The meaning of these equivalences is explained in the following. They display as:
\begin{align}{}
f&\equiv_{\eta}x\mapsto f(x)\label{eq:lambda-eta}\\
x\mapsto f(x)&\equiv_{\alpha}y\mapsto f(y)\label{eq:lambda-alpha}\\
\left(y\mapsto \mathbf{term}\right)(x)&\equiv_{\beta}\mathbf{term}[y:=x]\label{eq:lambda-beta}\\
f(g(x))&\equiv_{\circ}(f\circ g)(x).\label{eq:lambda-comp}
\end{align}
The intention of stating these rules is to be able to distinguish and name them. Our application of the $\eta$-equivalence on univariate functions (\ref{eq:lambda-eta}) states, that a function $f$ and the $\lambda$-abstraction\footnote{In programming languages this is sometimes called a \textit{wrapper function}.} immediately applying the argument $x\mapsto f(x)$ are computationally equivalent and therefore can be substituted against each other respecting the computation's result. The $\alpha$-equivalence (\ref{eq:lambda-alpha}) in this case states, that it does not matter for the computation how the argument is named, of course. So every time $\equiv_{\alpha}$ appears, the left hand side can be transformed in a computationally equivalent way to the right hand side by argument-renaming and vice versa. The $\beta$-equivalence (\ref{eq:lambda-beta}) expresses that an application of function $y\mapsto\mathbf{term}$, i.e. the $\mathbf{term}$ regarded as dependent on its variable $y$, to the argument $x$ is computationally equivalent to a $\mathbf{term}[y:=x]$ where all occurrences of $y$ are substituted for $x$. This is denoted by the substitution $[y:=x]$ acting on the $\mathbf{term}$ as a postfix operation. Lastly, not that much a rule of $\lambda$-calculus but more a definition of the composition operation $\circ$, is the rule (\ref{eq:lambda-comp}).

These rules (\ref{eq:lambda-eta}-\ref{eq:lambda-comp}) are somewhat standard rules that are most likely fulfilled in any context of computation. In $\lambda$-calculus every function takes exactly one argument and has one result which is a perfect interpretation for \textit{univariate} calculus. In computational electromagnetism, the representations of the considered objects, the electric and magnetic fields, the geometry, e.g. when given by parametrized coordinates, and coordinate transformations are expressed as \textit{multivariate} functions, taking multiple arguments to multiple results\footnote{We use the nomer multivariate, although it usually denotes functions taking multiple arguments to one result. Since in our case the results are not correlated to each other, and functions that give multiple uncorrelated results can be represented as a collection of these multivariate functions in the usual sense, we do not distinguish the terms here that much.}. Multiple arguments can be already thought of being represented as one argument with the help of the notion of a tuple, where the single arguments are separated by commas. Multiple results can be thought of as tuples in a similar manner. Yet, we choose a notation here that allows a multiple-argument-interpretation instead of tuples. It seems most familiar to the engineering community and does not pose a limitation since a multiple-argument-interpretation is translatable to a one-argument-interpretation.

That notation is motivated by the tediousness of multivariate calculus to express function application for these multiple arguments\footnote{Our proposed variant is mostly borrowed from the \textit{parameter-pack expansion} which is a carefully specified notation that appeared in the standard of the C++ programming language\cite{cpp17} first in its 2011 version. A parameter-pack can only appear in a meta-computation expressed within the templating system of C++. This notation is implemented in all current compilers complying to that standard.}. For a $\mathbf{term}$ we denote the \textit{expansion} by $\mathbf{term}...$ which should be computationally equivalent to a context where the comma-separation of copies of the term substituted with every single parameter, or variable in our case, of a tuple is applied. If $\underline{x}$ denotes a tuple of four parameters, the expansion of the most simple term, just consisting of $\underline{x}$ itself, corresponds to
\begin{equation*}
(\underline{x}...)\equiv_{m} (x^1, x^2, x^3, x^4)\label{eq:ppexpansion}\;.
\end{equation*}
Here, the tuple expansion $...$ captures\footnote{Just like a quantifier the tuple expansion binds unbound tuples where the unbound tuples are underlined.} all tuples in the term $\underline{x}$, which is just $\underline{x}$, and expands them to $x^1, x^2, x^3, x^4$ within the original term to produce the resulting term.
The three dots are used frequently in a meta-logical manner where it is clear from the context how to continue the pattern. When it comes to an implementation, one needs to make this pattern-repetition precise. In the following we make use of the three dots $...$ only in the sense of this kind of expansion, where the tuple is again underlined to highlight its meaning as a placeholder. The unexpanded term is denoted in an $m$-way  as computationally equivalent $\equiv_{m}$ to the expanded one.

The reason for introducing this particular notation is that it supports us in making precise arguments about multivariate functions in the previous sense. Our most important application is to express multivariate function application. E.g. suppose $g$ is a multivariate function in $\mathbb{R}^2\rightarrow\mathbb{R}^3$ such that it can be decomposed into functions $g^1$, $g^2$ and $g^3$ in $\mathbb{R}^2\rightarrow\mathbb{R}$, then we have two computationally equivalent terms with the \textit{nested use of the operation of tuple-expansion} $...$ :
\begin{align*}
\left({\color{green}\underline{g}}\left({\color{blue}\underline{x}...}\right){\color{green}...}\right)
&\equiv_{m}\left({\color{green}\underline{g}}\left({\color{blue}x^1},{\color{blue}x^2}\right){\color{green}...}\right)\\
&\equiv_{m}\left({\color{green}g^1}\left({\color{blue}x^1},{\color{blue}x^2}\right),{\color{green}g^2}\left({\color{blue}x^1},{\color{blue}x^2}\right),{\color{green}g^3}\left({\color{blue}x^1},{\color{blue}x^2}\right)\right).
\end{align*}
Here, the green tuple expansion ${\color{green}...}$ captures the green tuple ${\color{green}\underline{g}}$ where the blue tuple expansion ${\color{blue}...}$ captures the blue tuple ${\color{blue}\underline{x}}$.
Another use\footnote{This use is also borrowed from the C++ programming language standard.} is, given that $\gamma$ is a multivariate function in $\mathbb{R}^1\rightarrow\mathbb{R}^3$ that can be decomposed into the functions $\gamma^1$, $\gamma^2$ and $\gamma^3$ in $\mathbb{R}^1\rightarrow\mathbb{R}$, then we have two computationally equivalent terms with the expansion $...$ of \textit{multiple nested tuples} $\underline{\gamma}$ and $\underline{x}$:
$$\left({\color{blue}\underline{\gamma}}\left({\color{blue}\underline{x}}\right){\color{blue}...}\right)\equiv_{m} \left({\color{blue}\gamma^1}\left({\color{blue}x^1}\right),{\color{blue}\gamma^2}\left({\color{blue}x^2}\right),{\color{blue}\gamma^3}\left({\color{blue}x^3}\right)\right).$$
Here, the blue tuple expansion ${\color{blue}...}$ captures both blue tuples ${\color{blue}\underline{\gamma}}$ and ${\color{blue}\underline{x}}$.
Given the notion of tuples $\underline{x}$, $\underline{y}$ and the operation of tuple-expansion $...$, we can restate the previous computational equivalences (\ref{eq:lambda-eta}-\ref{eq:lambda-comp}) in their multivariate version (\ref{eq:lambda-eta-m}-\ref{eq:lambda-comp-m}):
\begin{align}{}
f&\equiv_{\eta}(\underline{x}...)\mapsto f(\underline{x}...)\label{eq:lambda-eta-m}\\
(\underline{x}...)\mapsto f(\underline{x}...)&\equiv_{\alpha}(\underline{y}...)\mapsto f(\underline{y}...)\label{eq:lambda-alpha-m}\\
\left((\underline{y}...)\mapsto \mathbf{term}\right)(\underline{x}...)&\equiv_{\beta}\mathbf{term}\left[(\underline{y}:=\underline{x})...\right]\label{eq:lambda-beta-m}\\
f(\underline{g}(\underline{x}...)...)&\equiv_{\circ}(f\circ g)(\underline{x}...)\label{eq:lambda-comp-m}\;.
\end{align}
Note especially how expansion interacts with composition of multivariate functions in (\ref{eq:lambda-comp-m}).

One more remark about tuples: you might have noticed that, despite underline and dots, the rules (\ref{eq:lambda-eta-m}-\ref{eq:lambda-comp-m}) exactly match the rules (\ref{eq:lambda-eta}-\ref{eq:lambda-comp}). That is not a coincidence.
Indeed, we could identify a scalar, with the one-tuple of scalars and have just one \textit{generalized} version of the rules. This works for all tuples, including one-tuples, and therefore also for all scalars.
Also, in Sec. \ref{sec:theory} we do regard multivariate functions as mapping tuples of scalars to tuples again\footnote{This means especially that we do not make use of \textit{currying} to express the multivariate functions. Furthermore as it turns out, the necessary (tuple-) arity of functions within this paper is always one.}.
Different notation and a reference to parameter-packs are just given, to support an interpretation within a language that does not identify scalars with one-tuples and may even distinguish tuples from a list of function-arguments\footnote{which is usually the case in programming languages tagged \textit{imperative}}.

You are free to ignore the three dots, when targeting an interpretation where function arguments and tuples are treated the same way\footnote{which is usually the case in programming languages tagged \textit{functional}: a function takes exactly one argument, which might be a tuple, and there is no difference of \textit{something} and the one-tuple of \textit{something}}. But as with higher order functions, tuples add a small burden on the learning curve and it is sometimes convenient to just think of a written-out version when comparing with the literature. One can test this preference by looking at $f\,(x,y,z)$ and if that should be \textit{a function $f$, applied to the function arguments $x$, $y$ and $z$}, then the tuple-expansion notation might be a fit. But if you are comfortable with $(x,y,z)$ being a tuple and if that tuple would be named $\tau$ your preferred notation is just $f\,\tau$ then you might want to ignore the dots. When programming, this choice is made by the programming language.

\subsection{Encoding the partial derivative}

We make use of the previously introduced equivalences to formulate what a partial derivative should be in that context. It is thought of as being the univariate derivative of a multivariate function which is regarded as a univariate function only depending on its one argument that we are taking the derivative of. That \textit{univariate regarding of a multivariate function} can be made precise now:
\begin{align}
h&\equiv_{\eta}(\underline{x}...)\mapsto h(\underline{x}...)\label{eq:c1}\\
&\equiv_{\beta}(\underline{x}...)\mapsto \left(x^2 \mapsto h(\underline{x}...)\right)(x^2)\label{eq:c2}\\
&\equiv_{m}({\color{blue}x^1},{\color{blue}x^2},{\color{blue}x^3})\mapsto \left({\color{purple}x^2} \mapsto h({\color{blue}x^1},{\color{purple}x^2},{\color{blue}x^3})\right)({\color{blue}x^2})\label{eq:c3}\\
&\equiv_{\alpha}({\color{blue}x^1},{\color{blue}x^2},{\color{blue}x^3})\mapsto \left({\color{purple}\mathbf{z}} \mapsto h({\color{blue}x^1},{\color{purple}\mathbf{z}},{\color{blue}x^3})\right)({\color{blue}x^2})\label{eq:c4}\\
&\equiv_{m}({\color{blue}\underline{x}}...)\mapsto \left({\color{purple}\mathbf{z}} \mapsto h\left({\color{blue}\underline{x}}...\left[\bullet^2:={\color{purple}\mathbf{z}}\right]\right)\right)({\color{blue}x^2}).\label{eq:c5}
\end{align}
Suppose the multivariate function $h$ is in $\mathbb{R}^3\rightarrow\mathbb{R}$. Then $h$ is computationally equivalent in an $\eta$-way to the multivariate function $(\underline{x}...)~\mapsto~ h(\underline{x}...)$ as in (\ref{eq:c1}). Just the inner $\mathbf{term}$ $h(\underline{x}...)$ of that new multivariate function is computationally equivalent to $\left(x^2 \mapsto h(\underline{x}...)\right)(x^2)$ in a univariate-$\beta$-way (\ref{eq:c2}). To see this, for the example, we look at the expanded version (\ref{eq:c3}). What happened is that the inner abstraction of ${\color{purple}x^2}$ is \textit{shadowing}\footnote{In theoretical computer science this is usually realized not by shadowing, but by limiting the $\alpha$-equivalence to the cases where the argument $x$ of $x\mapsto\mathbf{term}$ does not occur as a free variable of the $\mathbf{term}$, which is stated as $x\notin \text{FV}(\mathbf{term})$. But shadowing exists in the most programming languages.} the outer argument ${\color{blue}x^2}$. To highlight this difference, we explicitly rename the inner ${\color{purple}x^2}$ into an $\alpha$-equivalent function with ${\color{purple}\mathbf{z}}$ occuring instead (\ref{eq:c4}). This in a multivariate way constitutes the \textit{substituted expansion of the tuple} ${\color{blue}\underline{x}}$, denoted as ${\color{blue}\underline{x}}...\left[\bullet^2:={\color{purple}\mathbf{z}}\right]$, where entry 2 is replaced with ${\color{purple}\mathbf{z}}$ as in (\ref{eq:c5}).

\begin{figure}[h]
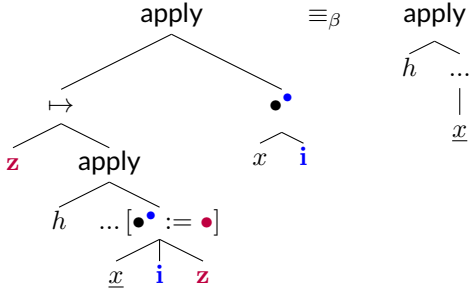
\centering
	\Tree[.apply [.$\mapsto$ [.${\color{purple}\mathbf{z}}$ ] [.apply [.$h$ ] [.$...\left[\bullet^{{\color{blue}\bullet}}:={\color{purple}\bullet}\right]$ [.$\underline{x}$ ] [.${\color{blue}\mathbf{i}}$ ] [.${\color{purple}\mathbf{z}}$ ] ]] ] [.$\bullet^{\color{blue}\bullet}$ $x$ ${\color{blue}\mathbf{i}}$ ]] $\equiv_{\beta}$ \Tree[.apply [.$h$ ] [.$...$ [.$\underline{x}$ ]] ]
	\caption{Univariate abstraction over an \mbox{expansion $...$} is expressed explicitly as substituted expansion $...\left[\bullet^{{\color{blue}\bullet}}:={\color{purple}\bullet}\right]$ to avoid implicit shadowing.}
	\label{fig:syntax-tree}
\end{figure}

This leads to the last rule of computational equivalence that we need for our considerations and it relates a multivariate function application to the use of a univariate function application:
\begin{equation}
\left({\color{purple}\mathbf{z}}\mapsto h\left(\underline{x}...\left[\bullet^{{\color{blue}\mathbf{i}}}:={\color{purple}\mathbf{z}}\right]\right)\right)\left({x}^{{\color{blue}\mathbf{i}}}\right)\equiv_{\beta}h\left(\underline{x}...\right)\label{eq:multi-uni}.
\end{equation}
To better familiarize with it, looking forward to an implementation, we give the syntax tree of this rule in Fig. \ref{fig:syntax-tree}.

\hyphenation{la-ter-on}

That is, finally, enough to define the partial derivative on multivariate functions $f:\mathbb{R}^d\rightarrow\mathbb{R}^c$ by the notion of the derivative ${\color{red}\mathbf{'}}$ on univariate functions. For a general arity and the indices ${\color{blue}\mathbf{j}}\in[1,c]$ and ${\color{blue}\mathbf{i}}\in[1,d]$ it is given as the multivariate function
\begin{equation}
\frac{\partial f^{\color{blue}\mathbf{j}}}{\partial x^{\color{blue}\mathbf{i}}}\coloneqq
({\underline{x}}...)\mapsto \left({\color{purple}\mathbf{z}} \mapsto f^{\color{blue}\mathbf{j}}\left({\underline{x}}...\left[\bullet^{\color{blue}\mathbf{i}}:={\color{purple}\mathbf{z}}\right]\right)\right)^{\color{red}\mathbf{'}}({x^{\color{blue}\mathbf{i}}})\label{eq:partial-derivative}\;,
\end{equation}
\noindent where $f^{\color{blue}\mathbf{j}}$ is the projection $\text{proj}^{\color{blue}\mathbf{j}}\circ f$ of the ${\color{blue}\mathbf{j}}$-th result of the multivariate function $f$ or similiarly the ${\color{blue}\mathbf{j}}$-th part of the decomposition of $f$ in the previously discussed manner.

In that definition (\ref{eq:partial-derivative}) we do not use the information about how to name the argument, with respect to which we are taking the partial derivative.
That is the case because partial derivatives with differently named arguments are computationally equivalent by $\alpha$-equivalence:
$$\frac{\partial f^{\color{blue}\mathbf{j}}}{\partial x^{\color{blue}\mathbf{i}}}\;\equiv_{\alpha}\;\frac{\partial f^{\color{blue}\mathbf{j}}}{\partial y^{\color{blue}\mathbf{i}}}\;\equiv_{\alpha}\;\partial_{\color{blue}\mathbf{i}}f^{\color{blue}\mathbf{j}}\;.$$
There are two remarks here to make. Firstly, the computational equivalence of the partial derivative under renaming of the argument, i.e. the $\alpha$-equivalence, motivates to omit the variable name $\partial_{\color{blue}\mathbf{i}}f^{\color{blue}\mathbf{j}}$. Later-on, however, in the theory of differential forms, this exact spot to give a name to the argument is often used to indicate which charts are involved in the process of coordinate transition\footnote{One distinguishes the function-level partial derivative $\partial_i$ with respect to the $i$-th argument of a function from the vector field $\partial/\partial x^i$ induced by the $i$-th coordinate $x^i$, where both fulfill the rules of what it means to be called a \textit{derivative}.}.
That characteristic results from the use of function-abstraction to express the partial derivative instead of introducing a new \textit{form of expression} as in the example of P. Martin-L\"of given in Sec. \ref{sec2b}.
In this way, the definition (\ref{eq:partial-derivative}) does not bind any free variables of its argument-terms.

Secondly, for a transition along $f$ from $A$-coordinates to $B$-coordinates, i.e., where $f$ is a function expressing the $B$-coordinates in terms of $A$-coordinates \mbox{$(\underline{f}(\underline{a}...)...)=(\underline{b}(\underline{a}...)...)$}, we have $\frac{\partial f^{\color{blue}\mathbf{j}}}{\partial a^{\color{blue}\mathbf{i}}}(\underline{a}...)$ to constitute the number in the ${\color{blue}\mathbf{j}}$'th row and the ${\color{blue}\mathbf{i}}$'th column of the Jacobi-matrix $\mathbf{J}_f$ evaluated in $A$-coordinates at $(\underline{a}...)$. That matrix is used to transform the numbers $(\underline{v_B}...)$ that are the vector-components with respect to the $B$-induced basis at a point given by the same $A$-coordinates into the the numbers $(\underline{v_A}...)$ that are the vector-components with respect to the $A$-induced basis at the same physical point by matrix-vector-multiplication \footnote{This is just the other way around as for the \textit{basis}, where $J_f$ transforms the $A$-induced basis into the $B$-induced basis.}. This scrutiny forms the foundation of a matrix-translation in terms of the Jacobi-matrix for different kinds of vectors. It is important to gain any support from encoding this logic into the notation and into the program to handle these different calculations and check them for consistency.

\subsection{The Chain-Rule revised}\label{sec3}


Using just these established conditions, we will derive what it means to have a notion of a chain rule for the partial derivative, lifting the notion of the univariate chain rule to the multivariate level. The whole calculation is given in appendix A. In order to create the multivariate listing in appendix A and the corresponding one for a concrete two-variate case in appendix B, we have implemented the tuple expansion the previously introduced way.



We begin in (\ref{eq:A1}) with the partial derivative that can be represented in an implementation not carrying anymore information than written in (\ref{eq:A1}), i.e. which function $f^{\color{blue}\mathbf{j}}\circ g$ it applies to with respect to which entry~${\color{blue}\mathbf{i}}$ or directly as the function that we encoded definitionally in (\ref{eq:partial-derivative}). In the first case an implementation needs to provide a function that converts these bits of information into that encoding. In the second case we directly operate on these objects. The multivariate function (\ref{eq:A2}) again does not need more information encoded than written out there and the data structure is very similar to the one resulting from a tree-like encoding of figure \ref{fig:syntax-tree}. The expanded terms for the two-variate case where ${\color{blue}\mathbf{i}}=1$ is given by (\ref{eq:B2}) and you can follow the expanded variant in appendix B alongside this investigation.

An equivalent computation (\ref{eq:A3}) is given by the multivariate \mbox{$\circ$-equivalence}, applying $f^{\color{blue}\mathbf{j}}$ to $\underline{g}$ instead of composing it with $g$. At this point, we make use of a linearity-property which needs to be fulfilled for a concrete realization of the univariate derivative ${\color{red}{'}}$ later-on. Namely that the univariate derivative of a multiply occuring argument is given by the sum of the univariate derivatives of each occurrence. We denote this by $=_{\text{lin}'}$ for the two-variate example given by:
\begin{align}%
&{\left({\color{purple}\mathbf{z}}\mapsto h\left({\color{purple}\mathbf{z}},{\color{purple}\mathbf{z}}\right)\right)}^{\color{red}\mathbf{'}}\left({x}\right)\nonumber\\
=_{\text{lin}'}\quad&{\left({\color{purple}\mathbf{z}}\mapsto h\left({\color{purple}\mathbf{z}},x\right)\right)}^{\color{red}\mathbf{'}}\left({x}\right)\label{eq:lin}\\
+\;&{\left({\color{purple}\mathbf{z}}\mapsto h\left(x,{\color{purple}\mathbf{z}}\right)\right)}^{\color{red}\mathbf{'}}\left(x\right)\nonumber\;.
\end{align}
\pagebreak
\newpage
\noindent For our general multivariate notation, $h$ has to be identified with
$$h:=({\underline{z}...})\mapsto f\left(\underline{g}\left(\underline{x}...\left[\bullet^{{\color{blue}\mathbf{i}}}:={\underline{z}}\right]\right)...\right),$$
\noindent leading to the general multivariate variant of this linearity, expressed with a summation ${\color{green}\sum_k}$ over a new \mbox{index ${\color{green}k}$:}
\begin{align}
&{\left({\color{purple}\mathbf{z}}\mapsto f\left(\underline{g}\left(\underline{x}...\left[\bullet^{{\color{blue}\mathbf{i}}}:={\color{purple}\mathbf{z}}\right]\right)...\right)\right)}^{\color{red}\mathbf{'}}\left({x}^{{\color{blue}\mathbf{i}}}\right)\label{eq:lin-gen}\\
=_{\text{lin}'}{\color{green}\sum_k}&{\left({\color{purple}\mathbf{z}}\mapsto f\left(\underline{g}\left(\underline{x}...\right)...\left[\bullet^{{\color{green}k}} :=g^{\color{green}k}\left(\underline{x}...\left[\bullet^{{\color{blue}\mathbf{i}}}:={\color{purple}\mathbf{z}}\right]\right)\right]\right)\right)}^{\color{red}\mathbf{'}}\left({x}^{{\color{blue}\mathbf{i}}}\right)\kern-0.25em,\nonumber
\end{align}
\noindent which expands in the two-variate case for ${\color{blue}\mathbf{i}}=1$ to:
\begin{align*}
&{\left({\color{purple}\mathbf{z}}\mapsto f\left(g^1\left({\color{purple}\mathbf{z}},x^2\right),g^2\left({\color{purple}\mathbf{z}},x^2\right)\right)\right)}^{\color{red}\mathbf{'}}\left(x^1\right)\\
=_{\text{lin}'}\quad&{\left({\color{purple}\mathbf{z}}\mapsto f\left(g^1\left({\color{purple}\mathbf{z}},x^2\right),g^2\left(x^1,x^2\right)\right)\right)}^{\color{red}\mathbf{'}}\left(x^1\right)\\
+\;&{\left({\color{purple}\mathbf{z}}\mapsto f\left(g^1\left(x^1,x^2\right),g^2\left({\color{purple}\mathbf{z}},x^2\right)\right)\right)}^{\color{red}\mathbf{'}}\left(x^1\right).
\end{align*}
Note the nested substitution in the right-hand-side term of (\ref{eq:lin-gen}) now, where only the application of the ${\color{green}k}$'th decomposition of $g$ is differently applied to the $x$'s of which just the ${\color{blue}\mathbf{i}}$'th one is replaced with ${\color{purple}\mathbf{z}}$. Therefore the linearity $=_{\text{lin}'}$ justifies whether (\ref{eq:A4}) computes the same result.
\begin{align}
&{\left({\color{purple}\mathbf{z}}\mapsto f\left(\underline{g}\left(\underline{x}...\right)...\left[\bullet^{{\color{green}k}} :=g^{\color{green}k}\left(\underline{x}...\left[\bullet^{{\color{blue}\mathbf{i}}}:={\color{purple}\mathbf{z}}\right]\right)\right]\right)\right)}\label{eq:nested-subst}\\
\equiv_{\circ}\;&{\left({\color{purple}\mathbf{z}}\mapsto f\left(\underline{g}\left(\underline{x}...\right)...\left[\bullet^{{\color{green}k}} :={\color{purple}\mathbf{z}}\right]\right)\right)}\circ
\left({\color{purple}\mathbf{z}}\mapsto
g^{\color{green}k}\left(\underline{x}...\left[\bullet^{{\color{blue}\mathbf{i}}}:={\color{purple}\mathbf{z}}\right]\right)
\right)
\nonumber
\end{align}
The nested substitution is computationally equivalent to the composition of univariate functions containing just a single substitution as in (\ref{eq:nested-subst}) which is the needed transformation that leads to (\ref{eq:A5}).

At this point, we have encoded the sum of ${\color{green}k}$ different univariate derivatives of a composition of two univariate functions (\ref{eq:A6}), where ${\color{green}k}$-times the univariate chain rule can be applied (\ref{eq:A7}) to lead to (\ref{eq:A8}). For the right multiplicand after transforming it in a $\beta$-way to the computationally equivalent form in (\ref{eq:A9}) it matches the definition of the partial derivative on $g$ (\ref{eq:A10}). The left multiplicand can be turned in a $\circ$- and $\beta$-way to the computationally equivalent form (\ref{eq:A11}-~\ref{eq:A12}) where the definition of the partial derivative again applies. This leads to the common form (\ref{eq:A13}) of the right hand side of the chain rule for the partial derivative of the composition of two functions $f^{\color{blue}\mathbf{j}}\circ g$, almost, but not quite:
\begin{equation}
\frac{\partial \left(f^{\color{blue}\mathbf{j}}\circ g\right)}{\partial {x}^{{\color{blue}\mathbf{i}}}}
=
\left(\underline{x}...\right)\mapsto {\color{green}\sum_k}\frac{\partial f^{\color{blue}\mathbf{j}}}{\partial y^{\color{green}k}}\left(\underline{g}\left(\underline{x}...\right)...\right)\cdot \frac{\partial g^{\color{green}k}}{\partial {x}^{{\color{blue}\mathbf{i}}}}\left(\underline{x}...\right)\;.\label{eq:chainrule}
\end{equation}
The applied calculus enforced an explicit mentioning of the abstraction $(\underline{x}...)\mapsto$ since these are function objects and only if they are applied to the same arguments, the one resulting number is equal for both sides:
\begin{equation}
\frac{\partial \left(f^{\color{blue}\mathbf{j}}\circ g\right)}{\partial {x}^{{\color{blue}\mathbf{i}}}}(\underline{x}...)
=
{\color{green}\sum_k}\frac{\partial f^{\color{blue}\mathbf{j}}}{\partial y^{\color{green}k}}\left(\underline{g}\left(\underline{x}...\right)...\right)\cdot \frac{\partial g^{\color{green}k}}{\partial {x}^{{\color{blue}\mathbf{i}}}}\left(\underline{x}...\right)\label{eq:chainrule-local}\;.
\end{equation}

\subsection{Targeting Tensor Calculus}\label{sec6}


In this paper our focus is to establish the lower interface that an encoding of the chain rule of multivariate functions demands from an encoding of the univariate chain rule.
It was investigated, how to define the partial derivative in computational terms.
We have shown in Sec. \ref{sec3} that this computational context is capable of deriving a chain rule for this definition.
In Sec. \ref{sec:theory} we will introduce an augmented $\lambda$-calculus based on the requirement to express a derivation of the chain rule from Sec. \ref{sec3}.
What remains open for discussion is the question whether that augmented $\lambda$-calculus is suitable to express definitions and derivations from tensor calculus.
It is also not obvious how the upper interface to tensor calculus should look like.
This section motivates why we think that our approach is extendable to express derivations from tensor calculus.

Continuing on (\ref{eq:A13}), with the $\circ$-equivalence we have a context (\ref{eq:A14}) where it is possible to make use of a function-level multiplication $\otimes$ that is given by the corresponding point-wise multiplication (\ref{eq:A15}). This is a binary operation and could be precomposed with a function applying $g$ to the left argument and the identity $\text{id}$ to the right argument. Defining such function is in favor for having just one binary operation on the two partial derivatives (\ref{eq:A16}), making a corresponding data structure definition even more obvious. Establishing a function-level summation $\oplus$ makes it possible to express the chain-rule in a completely so-called \textit{point-free}\footnote{i.e. a style where no arguments $(\underline{x}...)$ are present} style (\ref{eq:A17}). The objects reasoned about in this expression should correspond (denoted by $\cong_T$) to objects of the expression (\ref{eq:A18}) of tensor calculus, where unfortunately ${}'$ is a decoration on indices and not to be confused with the univariate derivative. We think that based on the way of that correspondence $\cong_T$ the question of encoding could be answered in a tractable way.

\pagebreak
What are the objects of tensor calculus that are common to reason about in computational electromagnetism? In the appendix of his book\cite{Tonti2013}, E. Tonti collects the notions of:
\begin{itemize}
	\item tensors and pseudotensors, such as tensor densities and tensor capacities, that differ in their transformation laws on a power of the determinant of the coordinate transition function,
	\item natural, reciprocal and physical basis vectors, leading to contravariant, covariant and physical components that are number-representations of various kinds of scalars and vectors in electromagnetic theory, and
	\item algebraic and metric dual vectors that constitute different representations of antisymmetric tensors.
\end{itemize}
\begin{figure}[h]
	\begin{align*}
	\text{Tonti}\;\text{2013}\;&{x'}^{h}=f^h(x^k) & \text{this}&\;\text{paper}\; f:A\rightarrow B \\
	\lambda^k_h(x')\stackrel{\text{def}}{=}\;
	&\frac{\partial x^k(x')}{\partial {x'}^h}
	&\text{J}^{{{\color{blue}\mathbf{k}}}}_{{{\color{blue}\mathbf{h}}'}}(\underline{a}...)
	&\cong_{T}\frac{\partial f^{\color{blue}\mathbf{k}}}{\partial {a}^{\color{blue}\mathbf{h}}}(\underline{a}...)
	\\
	\Lambda^h_k(x)\stackrel{\text{def}}{=}\;
	&\frac{\partial {x'}^h(x)}{\partial {x}^k}
	&\text{J}^{{{\color{blue}\mathbf{h}}'}}_{{{\color{blue}\mathbf{k}}}}(\underline{b}...)
	&\cong_{T}\frac{\partial {\left(f^{-1}\right)}^{\color{blue}\mathbf{h}}}{\partial {b}^{\color{blue}\mathbf{k}}}(\underline{b}...)
	\\
	\mathit{\Delta}(x')\stackrel{\text{def}}{=}\;&\text{det}(\lambda^k_h(x')) & (\textit{prop}&\textit{oses}\;\textit{permutations}) \\
	g\stackrel{\text{def}}{=}\;&\text{det}(g_{hk}) & (\textit{prop}&\textit{oses}\;\textit{permutations})
	\end{align*}%
	\caption{Denotational correspondences, where $g_{hk}$ is the metric tensor and $\text{det}$ being the determinant.}
	\label{fig:tonti}
\end{figure}
In classical electrodynamics, the physical base is often chosen because of its property to preserve the calculation for the length of a vector. This gives a direct interpretation for the measurement of such a quantity in a cartesian system, which is very valuable in a physical interpretation. These choices are combined with constructions such as the magnetic flux \textit{tuple} of numbers corresponding to the \textit{three-number representation} of the magnetic flux bi-covector at a point and similar constructions. Therefore, we think that it becomes arguable to investigate the computational aspects of such a correspondence. In accordance to follow his notation, which is very well chosen to support the application in various physical theories, we give correspondences in Fig. \ref{fig:tonti}.

Note especially, the choice of different symbols $\lambda$ and $\Lambda$ to reflect the information in which \textit{logical direction}\footnote{The direction, i.e. \textit{from the $A$ coordinate system to the $B$ coordinate system} or \textit{in the direction that $f$ is defined}, is meant here. To emphasize its distinction from the physical direction in space, we call it the \textit{logical} direction instead.} the partial derivative has to be taken and the drive to name the argument, $x$ or $x'$ respectively, to remember the coordinate transition function's domain. The difference between tensor calculus and the presented formalism is that we regard objects that are functions and function compositions where tensor calculus has a notion of \textit{coordinate system}. That is the key abstraction necessary to use in an implementation suitable of computing the chain-rule as a supporting layer. Consequently, we had no need to name the arguments and it is indeed not possible by $\alpha$-equivalence $\equiv_{\alpha}$ to encode that additional information.

Just to oppose it, we give in Fig. \ref{fig:tonti} another popular choice for denoting the partial derivative in tensor calculus  $\text{J}^{{{\color{blue}\mathbf{k}}}}_{{{\color{blue}\mathbf{h}}'}}$ for $\lambda^k_h$ and $\text{J}^{{{\color{blue}\mathbf{h}}'}}_{{{\color{blue}\mathbf{k}}}}$ for $\Lambda^h_k$. As mentioned before, the ${}'$ here should not be confused with the univariate derivative. The ${}'$ is a decoration on the indices ${\color{blue}\mathbf{k}}$ and ${\color{blue}\mathbf{h}}$ to represent their \textit{coordinate system belongingness}.
\begin{figure}
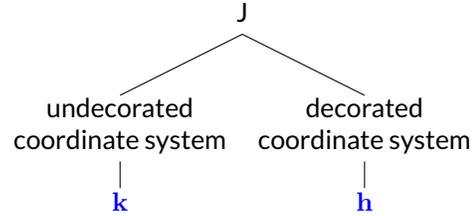
\centering%
	\Tree[.$\text{J}$ [.{undecorated\\coordinate system} ${\color{blue}\mathbf{k}}$ ] [.{decorated\\coordinate system} ${\color{blue}\mathbf{h}}$ ]]
	\caption{Encoding of the partial derivative used in tensor calculus}
	\label{fig:enc-tensor}
\end{figure}
Choosing different kinds of decorations for the indices to omit giving indices to the indices is an inevitable problem when multiple coordinate systems are considered.
In addition to that choice, there is the legitimate choice of the property of coordinate system belongingness being one of the index or being a property of the \textit{partial derivative object} itself.
The former perspective is taken in the notation we opposed which where the latter was denoted $\lambda$ or $\Lambda$ respectively.
This state of affairs is also shown in fig. \ref{fig:enc-tensor}. An answer to that question of choice highly influences the encoding of tensor calculus expressions for the purpose of an implementation.

\hyphenation{cla-ri-fied}
\hyphenation{de-co-ra-ted}

As promised in the title, we will show here transformation laws for the magnetic flux $B$ and the electric field $E$, although the reason of this paper is not the result but the process of deriving these laws. For a clarified choice of $\cong_{T}$, which we did not yet made in this paper, suppose that $Z$, $A$ and $B$ are given by left decorated $z{}^\backprime$, undecorated $a$ and right decorated $b'$ coordinate systems. In this notation, for clarification, the coordinate system belongingness is redundantly encoded in the choice of the letter, as well as in the decoration of that letter. This amounts to the habit that in the calculus of multivariate functions, just different letters are used, where in tensor calculus only different decorations are used. Then, for the two transition functions $g:Z\rightarrow A$ and $f:A\rightarrow B$ the tensor calculus expression that relates the covariant components $B_{{\color{blue}\mathbf{i}}'{\color{blue}\mathbf{j}}'}$ of the bi-covector of the right decorated coordinate system $b'$ to the ones $B_{{\color{green}\mathbf{i}}{\color{green}\mathbf{j}}}$ of the undecorated coordinate \mbox{system $a$} is given by:
\begin{equation*}
B_{{\color{blue}\mathbf{i}}'{\color{blue}\mathbf{j}}'}=\text{J}^{{\color{green}\mathbf{i}}}_{{\color{blue}\mathbf{i}}'}\text{J}^{{\color{green}\mathbf{j}}}_{{\color{blue}\mathbf{j}}'}B_{{\color{green}\mathbf{i}}{\color{green}\mathbf{j}}}\;,
\end{equation*}
where \textit{free} indices are highlighted in blue and \textit{bound} indices, which are summed over, are highlighted in green. This translates into:
\begin{equation*}
\begin{matrix*}[l]
{\displaystyle B_{{\color{blue}\mathbf{i}}'{\color{blue}\mathbf{j}}'}\equiv
\left(\underline{b}...\right)\mapsto
{\color{green}\sum_{{\color{green}\mathbf{i}}{\color{black}{}}}}
{\color{green}\sum_{{\color{green}\mathbf{j}}{\color{black}{}}}}}
&{\displaystyle
\frac{\partial \left(f^{-1}\right)^{{\color{green}\mathbf{i}}}}{\partial b^{{\color{blue}\mathbf{i}}'}}
\left(\underline{b}...\right)}\\
&{\displaystyle\cdot
\frac{\partial \left(f^{-1}\right)^{{\color{green}\mathbf{j}}}}{\partial {b}^{{\color{blue}\mathbf{j}}'}}\left(\underline{b}...\right)}\\[1em]
&{\displaystyle\cdot
B_{\color{green}\mathbf{i}\mathbf{j}}\left(\underline{f}(\underline{b}...)...\right)\;.}
\end{matrix*}
\end{equation*}
As $B_{\color{green}\mathbf{i}\mathbf{j}}$ should be regarded to \textit{naturally} live on the undecorated coordinates $a$ and the resulting object $B_{{\color{blue}\mathbf{i}}'{\color{blue}\mathbf{j}}'}$ to live on the right decorated coordinates $b'$, a precomposition with $f$ is necessary to obtain the $B_{\color{green}\mathbf{i}\mathbf{j}}$ value at $b'$ coordinates. Although this transformation goes in the same logical direction as the functions $g$ and $f$ are defined, the partial derivatives of inverses of these functions appear due to the contravariant transformation property of the considered electromagnetic quantity.

Tensor calculus is concerned about the invariance properties of different quantities. Suppose that two Jacobians cancel each other out in the following way
$$J^{{\color{blue}\mathbf{i}'}}_{{\color{green}\mathbf{i}}}J^{{\color{green}\mathbf i}}_{{\color{blue}\mathbf j'}}=\delta^{{\color{blue}\mathbf i'}}_{{\color{blue}\mathbf j'}}$$
where $\delta^{{\color{blue}\mathbf i'}}_{{\color{blue}\mathbf j'}}$ is the Kronecker delta which is $1$ for ${\color{blue}\mathbf i'}={\color{blue}\mathbf j'}$ and $0$ otherwise.
Then, it is \textit{easy} to see that the transformations of the tensor components in $S^{{\color{green}\mathbf{ij}}}T_{{\color{green}\mathbf i}}U_{{\color{green}\mathbf j}}$ will cancel each other out. If the independent transformations of $S^{{\color{blue}\mathbf{ij}}}$, $T_{{\color{blue}\mathbf i}}$ and $U_{{\color{blue}\mathbf j}}$ are
$$S^{{\color{blue}\mathbf{ij}}}=S^{{\color{green}\mathbf{i'j'}}}J^{{\color{blue}\mathbf i}}_{{\color{green}\mathbf i'}}J^{{\color{blue}\mathbf j}}_{{\color{green}\mathbf j'}},\quad\quad T_{{\color{blue}\mathbf i}}=T_{{\color{green}\mathbf k'}}J^{{\color{green}\mathbf k'}}_{{\color{blue}\mathbf i}}\quad\text{and}\quad U_{{\color{blue}\mathbf j}}=U_{{\color{green}\mathbf l'}}J^{{\color{green}\mathbf l'}}_{{\color{blue}\mathbf j}},$$
then we have for the composed term
\begin{align*}
S^{{\color{green}\mathbf{ij}}}T_{{\color{green}\mathbf i}}U_{{\color{green}\mathbf j}}&=\left(S^{{\color{green}\mathbf{i'j'}}}J^{{\color{green}\mathbf i}}_{{ \color{green}\mathbf i'}}J^{{\color{green}\mathbf j}}_{{\color{green}\mathbf j'}}\right)\left(T_{{\color{green}\mathbf k'}}J^{{\color{green}\mathbf k'}}_{{\color{green}\mathbf i}}\right)\left(U_{{\color{green}\mathbf l'}}J^{{\color{green}\mathbf l'}}_{{\color{green}\mathbf j}}\right)\\
&=S^{{\color{green}\mathbf{i'j'}}}\left(J^{{\color{green}\mathbf k'}}_{{\color{green}\mathbf i}}J^{{\color{green}\mathbf i}}_{{\color{green}\mathbf i'}}\right)\left(J^{{\color{green}\mathbf l'}}_{{\color{green}\mathbf j}}J^{{\color{green}\mathbf j}}_{{\color{green}\mathbf j'}}\right)T_{{\color{green}\mathbf k'}}U_{{\color{green}\mathbf l'}}\\
&=S^{{\color{green}\mathbf{i'j'}}}\left(\delta^{{\color{green}\mathbf k'}}_{{\color{green}\mathbf i'}}\right)\left(\delta^{{\color{green}\mathbf l'}}_{{\color{green}\mathbf j'}}\right)T_{{\color{green}\mathbf k'}}U_{{\color{green}\mathbf l'}}\\
&=S^{{\color{green}\mathbf{i'j'}}}T_{{\color{green}\mathbf i'}}U_{{\color{green}\mathbf j'}}.
\end{align*}
Here, the parentheses are present only for clarification since this tensor expression represents scalar multiplications.

In our previous example the transformation of a bi-covector needs to be \textit{justified} by invariance properties that are expressable within tensor calculus.
The appearance of Jacobians is due to this invariance. With our proposed formalism, we can express this as a computational equivalence.
If we were to have a representation for $B_{\color{blue}\mathbf{i}'\mathbf{j}'}$ and $B_{\color{green}\mathbf{i}\mathbf{j}}$ then such an equivalence should be derivable at this tensor-calculus-layer.
This might motivate that the statements that need to be proven, which arise in tensor calculus, are expressable in our proposed formalism.


Another example is the tensor calculus expression relating the contravariant components $E^{{\color{blue}\mathbf{i}}{}'}$ of vector $E$ in the right decorated coordinate system $b{}'$ to the ones $E^{{\color{green}\mathbf{i}}^\backprime}$ of the left decorated coordinate system $z^\backprime$ is given by:
\begin{equation*}
E^{{\color{blue}\mathbf{i}}'}
=\text{J}_{{\color{green}\mathbf{i}}}^{{\color{blue}\mathbf{i}}'}\text{J}_{{\color{green}\mathbf{i}}^\backprime}^{{\color{green}\mathbf{i}}}E^{{\color{green}\mathbf{i}}^\backprime}
=\text{J}_{{\color{green}\mathbf{i}}{}^\backprime}^{{\color{blue}\mathbf{i}}'}E^{{\color{green}\mathbf{i}}^{\backprime}}\;,
\end{equation*}
\noindent which translates into:
\begin{equation*}
\begin{matrix*}[l]
{\displaystyle
	E^{{\color{blue}\mathbf{i}}'}\equiv
	\left(\underline{b}...\right)\mapsto
	{\color{green}\sum_{{\color{green}\mathbf{i}}{\color{black}{}}}}
	{\color{green}\sum_{{\color{green}\mathbf{i}}{\color{black}{}^\backprime}}}
}&
{\displaystyle
	\;\frac{\partial f^{{\color{blue}\mathbf{i}}'}}{\partial a^{{\color{green}\mathbf{i}}}}
	\left(\underline{f}^{-1}\left(\underline{b}...\right)...\right)
}\\
&  \cdot
{\displaystyle
	\frac{\partial g^{{\color{green}\mathbf{i}}}}{\partial {z}^{{\color{green}\mathbf{i}}^\backprime}}
	\left(\underline{g}^{-1}(\underline{f}^{-1}(\underline{b}...)...)...\right)
}\\[1em]
&  \cdot
{\displaystyle
	E^{{\color{green}\mathbf{i}}^\backprime}
	\left(\underline{g}^{-1}(\underline{f}^{-1}(\underline{b}...)...)...\right)
}.
\end{matrix*}
\end{equation*}
\begin{equation*}
\begin{matrix*}[l]
{\displaystyle
	=_{\text{chain rule}}
	\left(\underline{b}...\right)\mapsto
	{\color{green}\sum_{{\color{green}\mathbf{i}}{\color{black}{}^\backprime}}}
}&
{\displaystyle
	\;\frac{\partial (f^{{\color{blue}\mathbf{i}}'}\circ g)}{\partial {z}^{{\color{green}\mathbf{i}}^\backprime}}
	\left(\underline{g}^{-1}(\underline{f}^{-1}(\underline{b}...)...)...\right)
}\\
&  \cdot
{\displaystyle
	E^{{\color{green}\mathbf{i}}^\backprime}
	\left(\underline{g}^{-1}(\underline{f}^{-1}(\underline{b}...)...)...\right)
}
\end{matrix*}
\end{equation*}
Here again the resulting $E^{{\color{blue}\mathbf{i}}'}$ should live on the right decorated coordinates $b'$, where the original $E^{{\color{green}\mathbf{i}}^\backprime}$ lives on the left decorated coordinates $z^\backprime$. The transformation happened again in the same logical direction as the functions go, but this time we have transformed twice. To apply the introduced partial differential of the multivariate functions, it becomes necessary to precompose with proper inverses to obtain an expression that again depends on the right decorated coordinates $b'$.

This second example of a chained coordinate transformation makes use of the derived chain rule which is justified in the current formalism. This should make it easy in an on-top tensor calculus layer to computationally proof 
$$
\text{J}_{{\color{green}\mathbf{i}}}^{{\color{blue}\mathbf{i}}'}\text{J}_{{\color{blue}\mathbf{i}}^\backprime}^{{\color{green}\mathbf{i}}}
=\text{J}_{{\color{blue}\mathbf{i}}{}^{\backprime}}^{{\color{blue}\mathbf{i}}'}
$$
when a clarified choice is made how the tensor calculus terms should correspond $\cong_{T}$ to terms from $\lambda$-calculus.

\section{Theoretical Framework}\label{sec:theory}

\subsection{Rules of inference}

A logical inference, e.g.,
\begin{prooftree}
\AxiomC{$A$}
\AxiomC{$B$}
\inference[2]{$A\;\&\; B$}
\end{prooftree}
\textit{``does not take us from the propositions $A$ and $B$ to the proposition $A\;\&\; B$''}\cite{martinlof83}. 
\textit{``Rather, it takes us from the affirmation of A and the affirmation of B to the affirmation of $A\;\&\; B$''}\cite{martinlof83}.
Making this explicit in writing could\footnote{P. Martin-L\"of attributes it to B. Russell, translating Frege's \textit{Urteil} into \textit{assertion}, and calling the combination of Frege's judgment stroke "$|$" and content stroke "$-$" the assertion sign "$\vdash$".} be
\begin{prooftree}
\AxiomC{$\vdash A$}
\AxiomC{$\vdash B$}
\RightLabel{.}
\inference[2]{$\vdash A\;\&\; B$}
\end{prooftree}
there is a need to distinguish two kinds of entities:
%
\begin{itemize}
\item \textit{``the entities that the logical operations operate on, which we call propositions''}\cite{martinlof83}~\footnote{that use of the word \textit{proposition} is again attributed to B. Russell}, which are \textit{``affirmed in an affirmation and denied in a denial''}\cite{martinlof83}.
\item and \textit{``the things that the logical laws, by which I mean the rules of inference, operate on, which we normally call assertions''}\cite{martinlof83}, which are \textit{``those that we prove and that appear as premises and conclusion of a logical inference''}\cite{martinlof83}.
\end{itemize}
We are examining that topic at this point in the paper for two reasons:
One is in preparation of stating \textit{introduction rules} for an augmented $\lambda$-calculus.
The other is, because the word \textit{proposition} has a different meaning in logic than it has in most of mathematics.

A logicians issue with the mathematical wording would be, that \textit{``a theorem is sometimes called a proposition, sometimes a theorem''}\cite{martinlof83}. And thus \textit{``we have two words for the things that we prove, proposition and theorem''}\cite{martinlof83}.
Now, \textit{``what we prove, in particular, the premises and conclusion of a logical inference''}\cite{martinlof83} are not called propositions, but judgments or assertions.


There is one technicality here that one might not even notice. Strictly speaking, the word judgement, or assertion, is used in particular for the premises and conclusion of a logical inference where it usually means an affirmation or denial. Most of modern logic gets along with just affirmations.
A formula is not affirmed directly, but it has to be \textit{grasped} as a proposition and that proposition then can be affirmed.
When
\begin{prooftree}
\AxiomC{$A\text{ prop}$}
\AxiomC{$B\text{ prop}$}
\AxiomC{$A\text{ true}$}
\inference[3]{$A\;\vee\; B\text{ true}$}
\end{prooftree}
should be a \textit{rule of disjunction	introduction}, then grasping $A$ and $B$ as propositions, $A\text{ prop}$ and $B\text{ prop}$ do figure as premises for that rule although they are not an affirmation nor a denial.
P. Martin-L\"of extends a use of the word \textit{judgment} to include such new \textit{forms of judgment} which are not only affirmations or denials anymore.
Extending that usage allows us to denote premises and conclusion of an inference as \textit{judgments of some specific form}.
This wording is important, because for typed $\lambda$-calculus our goal is the derivation of judgments the form $\Gamma\vdash T : a$ which means $T$ has type $a$ in context $\Gamma$. These judgments appear as the premises and conclusion of \textit{type checking rules of inference}.
There are three basic introduction type checking rules of typed $\lambda$-calculus: introducing $\lambda$-abstraction, introducing function application and introducing variable usage.


\subsection{Typed $\lambda$-calculus}

In order to formalize the previously motivated application in Sec. \ref{sec:problem}, we spoke about \textit{univariate} and \textit{multivariate} functions, \textit{tuples} made of \textit{scalars} and \textit{indices} for various operations on tuples and multivariate functions.
These all make valid \textit{types} in our consideration and therefore we model a ``Type'' in our augmented $\lambda$-calculus to be introduced by the following introduction rules:
\begin{prooftree}
	\AxiomC{$k : \mathbb{N}$}
	\RightLabel{($\Index -$)}
	\inference[1]{$\Index k : \text{ Type}$}
\end{prooftree}
\begin{prooftree}
	\AxiomC{}
	\RightLabel{($\funII$)}
	\inference[1]{$\funII : \text{ Type}$}
\end{prooftree}
\begin{prooftree}
	\AxiomC{$m : \mathbb{N}$}
	\RightLabel{($\funMI -$)}
	\inference[1]{$\funMI m : \text{ Type}$}
\end{prooftree}
\begin{prooftree}
	\AxiomC{$n : \mathbb{N}$}
	\RightLabel{($\funIN -$)}
	\inference[1]{$\funIN n : \text{ Type}$}
\end{prooftree}
\begin{prooftree}
	\AxiomC{$m : \mathbb{N}$}
	\AxiomC{$n : \mathbb{N}$}
	\RightLabel{($\funMN -\; -$)}
	\inference[2]{$\funMN m\;n : \text{ Type}$}
\end{prooftree}
\begin{prooftree}
	\AxiomC{$k : \mathbb{N}$}
	\RightLabel{($\tuple -$)}
	\inference[1]{$\tuple k : \text{ Type}$}
\end{prooftree}
\begin{prooftree}
	\AxiomC{}
	\RightLabel{($\scalar$).}
	\inference[1]{$\scalar : \text{ Type}$}
\end{prooftree}
These mean, that there \textit{is} a type for univariate functions ``fun11'' and a type for scalars. 
For every natural number there is one type of indices, one type of functions taking $m$ arguments to a single output ``funM1'' and one type of functions taking a single input to $n$ outputs ``fun1N''.
For every two natural numbers $m$ and $n$ there is a function type taking $m$ inputs to $n$ outputs ``funMN''.
These are purely syntactical introduction rules that are named semantically but do not have their intended meaning yet.
But this ``Type'' serves as index set over which we will define \textit{the family of valid terms} meaning we regard the totality of terms partitioned by their ``Type''.

These rules, in a very exact sense, correspond to a datatype definition in the Agda\cite{norell2009} language\footnote{
The Agda language, on the one hand can be introduced as a functional programming language that, on the other hand, is powerful enough to express constructive mathematics.
Agda builds on top of a type theory, as introduced by P. Martin-L\"of.
It supports dependently typed pattern matching, using so-called \textit{Miller pattern unification}, with Σ-types, inductive datatypes and universe polymorphism.}.

\begin{figure}[h]
\begin{Shaded}
\begin{Highlighting}[]
	\KeywordTok{data}\NormalTok{ Type }\OtherTok{:} \DataTypeTok{Set} \KeywordTok{where}
	\NormalTok{  index  }\OtherTok{:}\NormalTok{ ℕ }\OtherTok{→}\NormalTok{ Type}
	\NormalTok{  fun11  }\OtherTok{:}\NormalTok{ Type}
	\NormalTok{  funM1  }\OtherTok{:}\NormalTok{ ℕ }\OtherTok{→}\NormalTok{ Type}
	\NormalTok{  fun1N  }\OtherTok{:}\NormalTok{ ℕ }\OtherTok{→}\NormalTok{ Type}
	\NormalTok{  funMN  }\OtherTok{:}\NormalTok{ ℕ }\OtherTok{→}\NormalTok{ ℕ }\OtherTok{→}\NormalTok{ Type}
	\NormalTok{  tuple  }\OtherTok{:}\NormalTok{ ℕ }\OtherTok{→}\NormalTok{ Type}
	\NormalTok{  scalar }\OtherTok{:}\NormalTok{ Type}
\end{Highlighting}
\end{Shaded}
\caption{Agda datatype of custom types to be regarded in an augmented $\lambda$-calculus.}
\label{fig:customtypes}
\end{figure}

We have, that for every "Type" that can be introduced by our stated introduction rules, there is exactly one element in the datatype that we have defined within the Agda language and vice versa.
This property makes it suitable to support our formalization as we go along. Here $\mathbb{N}$ is inductively defined in the usual way which is not much of interest here.
For the formalization we introduced also the totality ``Name'' of names where variables are chosen from, but this also a minor point.


Some of the introduction type checking rules of $\lambda$-calculus in general are $\lambda$-abstraction and function application.
Where the latter usually is denoted just by juxtaposition, without an explicit operator, we emphasize this by the use of $\llparenthesis$ and $\rrparenthesis$.

The first rule of $\lambda$-abstraction for $a$ and $b$ being types in our $\lambda$-calculus and $x$ being a name is written as
\begin{prooftree}
	\AxiomC{$\Gamma , (x , a) \vdash T : b$}
	\RightLabel{($\lambda - {\;.\;} -$).}
	\inference[1]{$\Gamma \vdash \lambda x {\;.\;} T : a\to b$}
\end{prooftree}
It takes us from the judgment that \textit{\guillemotleft~in a context consisting of first, $\Gamma$ and second, the variable $x$ being of type $a$, within that context the term $T$ is of type $b$ \guillemotright} to the judgment that \textit{\guillemotleft~in context $\Gamma$ the term $\lambda\,x\,.\,T$ is of type $a\to b$~\guillemotright}.

In our application in Sec. \ref{sec:problem} we only needed to abstract over scalars or tuples, so a much stricter rule can be used for a formalization. We have chosen four much simpler rules instead which fix the types to the four combinations of tuples and scalars. You can find them in the appendix \ref{sec:allrules}.
The reason to take this simplification is that a following interpretation in section \ref{sec:software} will become easier having the types fixed. This is possible because we are not targeting to formalize a general purpose programming language, but rather a very specific one just targeting the partial derivative.

A second rule introduces function application:
\begin{prooftree}
	\AxiomC{$\Gamma \vdash T : a\to b$}
	\AxiomC{$\Gamma \vdash U : a $}
	\RightLabel{($- {\;} \llparenthesis - \rrparenthesis$)}
	\inference[2]{$\Gamma \vdash T \llparenthesis\; U \;\rrparenthesis : b$}
\end{prooftree}
It takes us from the two judgments that \textit{\guillemotleft~in a context $\Gamma$ the term $T$ is of type $a\to b$ \guillemotright} and \textit{\guillemotleft~in the same context $\Gamma$ the term $U$ is of type $a$ \guillemotright} to the judgment \textit{\guillemotleft~in context $\Gamma$ again, $T\llparenthesis U\rrparenthesis$ is of type $b$ \guillemotright}.
In our formalization we chose to have four such introduction rules operating on the corresponding argument types. One for each type of function.

\subsection{De Bruijn indices}

There is one very basic key technique to work out for developing sane introduction rules that really respect the typing of variables with respect to some \textit{context} $\Gamma$.
It is noteworthy, that this technique is necessary to produce correctness guarantees from a programming language's type checker as motivated in Sec. \ref{sec:computationalcontext}.
Unfortunately it is only expressable in a dependently typed programming language. In other programming languages, the following rules reduce to a list data structure.

As introduction rules for a context we chose that there \textit{is} an empty context
\begin{prooftree}
	\AxiomC{}
	\RightLabel{($[]$)}
	\inference[1]{$[] : \text{ Context}$}
\end{prooftree}
and, when given a context $\Gamma$, we can form a new one $\Gamma , (x , a)$ for every name-type combination $(x , a)$.
\begin{prooftree}
	\AxiomC{$\Gamma : \text{ Context}$}
	\AxiomC{$x : \text{ Name}$}
	\AxiomC{$a : \text{ Type}$}
	\RightLabel{($- , (-,-)$)}
	\inference[3]{$\Gamma , (x , a) : \text{ Context}$}
\end{prooftree}
These rules make a context to \textit{a list of tuples containing a name and a type} in our consideration.
The de Bruijn indices\cite{debruijn72} that we are going to work out will be indices that are guaranteed by their type to really point to a specific name-type combination within such context.
One could even model a context as a list of just types and without the names. Within the Agda formalization we found it very expressive to have this redundant piece of information available.
Still, a variable is to be identified by its de Bruijn index and not by its name.

The first rule introduces a judgment that \textit{\guillemotleft~de Bruijn index zero is an element of the type of de Bruijn indices that show the first name-type combination of a context being in that context \guillemotright}.
\begin{prooftree}
	\AxiomC{$\Gamma : \text{ Context}$}
	\AxiomC{$x : \text{ Name}$}
	\AxiomC{$a : \text{ Type}$}
	\RightLabel{($\text{zero}$)}
	\inference[3]{$\text{zero} : (x , a) \in \Gamma , (x , a)$}
\end{prooftree}
The second rule introduces a judgment that \textit{\guillemotleft~an incremented de Bruijn index shows that a name-type combination is part of an appended context, given that it did so before \guillemotright}.
\begin{prooftree}
	\AxiomC{$i^b : (x , a) \in \Gamma$}
	\AxiomC{$y : \text{ Name}$}
	\AxiomC{$b : \text{ Type}$}
	\RightLabel{($\text{suc} -$)}
	\inference[3]{$\text{suc } i^b : (x , a) \in \Gamma , (y , b)$}
\end{prooftree}
With these rules it is possible to give meaning to the last standard introduction type check rule of $\lambda$-calculus.
It introduces the judgment that \textit{\guillemotleft~in context $\Gamma$ the variable $x$ is of type $a$, because of~~$\because$~~the de Bruijn index $i^b$ \guillemotright}.
\begin{prooftree}
	\AxiomC{$i^b : (x , a) \in \Gamma$}
	\RightLabel{($- \because - $)}
	\inference[1]{$\Gamma \vdash x \;\because\; i^b : a$}
\end{prooftree}
Where this matches closely our Agda formalization, the rule is sometimes written more intuitively as
\begin{prooftree}
	\AxiomC{$x : a \in \Gamma$}
	\RightLabel{(Var)}
	\inference[1]{$\Gamma \vdash x : a$}
\end{prooftree}
or even
\begin{prooftree}
	\AxiomC{$x \in \Gamma$}
	\RightLabel{(Var).}
	\inference[1]{$\Gamma \vdash x : \Gamma(x)$}
\end{prooftree}

\subsection{Augmenting $\lambda$-calculus}

By just \textit{using} the $\lambda$-calculus we got into the previous three rules and the use of de Bruijn indices even without any specifics from our application.
After paying that entry fee for which do not exist many alternatives, we can finally work-in our application specific operations which are: \textit{substitiution of the $i$-th component}, \textit{projecting out the $i$-th component}, \textit{composition of functions}, \textit{the univariate derivative} and \textit{scalar multiplication}.

\begin{prooftree}
	\AxiomC{$\Gamma \vdash T : a $}
	\AxiomC{$i : \mathbb{N} $}
	\AxiomC{$\Gamma \vdash U : b $}
	\RightLabel{($- [\bullet - := -]$)}
	\inference[3]{$\Gamma \vdash T [\bullet i := U] : a$}
\end{prooftree}

\begin{prooftree}
	\AxiomC{$\Gamma \vdash T : a $}
	\AxiomC{$i : \mathbb{N} $}
	\RightLabel{($- \;\hat{}\; -$)}
	\inference[2]{$\Gamma \vdash T \;\hat{}\; i : b$}
\end{prooftree}

\begin{prooftree}
	\AxiomC{$\Gamma \vdash T : b\to c $}
	\AxiomC{$\Gamma \vdash U : a\to b $}
	\RightLabel{($- \circ -$)}
	\inference[2]{$\Gamma \vdash T \circ U : a\to c$}
\end{prooftree}

\begin{prooftree}
	\AxiomC{$\Gamma \vdash T : r\to r$}
	\RightLabel{($- {\;}^{'}$)}
	\inference[1]{$\Gamma \vdash T {\;}^{'} : r\to r$}
\end{prooftree}

\begin{prooftree}
	\AxiomC{$\Gamma \vdash T : r$}
	\AxiomC{$\Gamma \vdash U : r$}
	\RightLabel{($- \cdot -$)}
	\inference[2]{$\Gamma \vdash T \cdot U : r$}
\end{prooftree}

For the previously introduced types in our specific $\lambda$-calculus, we introduced six rules for substitution, six rules for projection, four obvious rules for composition of functions and one rule for the univariate derivative as well as one rule for scalar multiplication.
You can find the rules in the appendix \ref{sec:allrules} and the Agda datatype of all well-formed $\lambda$-calculus terms is in Fig. \ref{fig:agdaterm}.

\begin{figure*}[t]
\begingroup
\fontsize{8pt}{10pt}\selectfont
\input{agdaTerm.tex}
\endgroup
\caption{The datatype "Term" of well formed terms for an \textit{augmented} $\lambda$-calculus suitable to express the partial derivative within the Agda programming language. $\text{Fin }k$ is the type of natural numbers less than $k$, sometimes denoted $\mathbb{N}_{k}$ or $\mathbb{N}_{<k}$.}
\label{fig:agdaterm}
\end{figure*}

%
%
%

\subsection{Chain of Justification}\label{sec4}

All the data structures and data transformations described in Sec. \ref{sec:problem}, represent computations for the partial derivative function. 
But even after translating them into an augmented $\lambda$-calculus, they are not yet more than the mere skeletons carrying around meta-data.
All the transformations we now implement on these $\lambda$-terms which should respect this, yet hypothetical, computation are just operations transforming that meta-data.

The resulting $\lambda$-terms can only be turned into a computation when a lower layer, i.e. an implementation providing the univariate derivative and a representation of functions, providing these computations, is present such that the terms can be interpreted, i.e. turned into a computation and executed.

There are just a few properties even possible to be proven without further assumption at this high level. 
We have made the distinction between a computational equivalence~$\equiv$ that is justified within our investigation by the computational equivalences of the $\lambda$-calculus and the propositional equality $=$ that is used when a property of the univariate derivative ${\color{red}\mathbf{'}}$, that operation we presupposed for our whole consideration, was made use of.

For the computational equivalences $\equiv$ there is some chance to express those in terms of $\alpha$-conversion, $\beta$-reduction and $\eta$-reduction.
But the provability of the equivalences denoted by $=$ depends on the underlying \textit{interpretation}.


Consistency of computational equivalence resulting from the presented transformations depends on a consistent implementation of the considered layer, of course, and precisely on a consistent implementation of these two equality-transformations of the lower layer. These two equality-transformations are in some sense \textit{dependencies} of our considered layer. The benefit is that the implementation of the considered layer can be verified in a way independently from a lower level application increasing the overall trust and decomposing monolithic software ventures into more modular ones.
Similar to the two assumptions $=_{\text{lin}}$ and $=_{\text{chain}}$ of the univariate derivative, it is possible to determine additional assumptions that are necessary in proofs of additional theorems.
In Sec. \ref{sec:software} we give guidance how these rather abstract assumptions become more concrete with a chosen interpretation for the augmented $\lambda$-calculus.

\section{Software Architecture}\label{sec:software}

\begin{figure*}[ht!]
	\centering
	\def\svgwidth{0.85\textwidth}
	\begingroup
	\fontsize{8pt}{10pt}\selectfont
	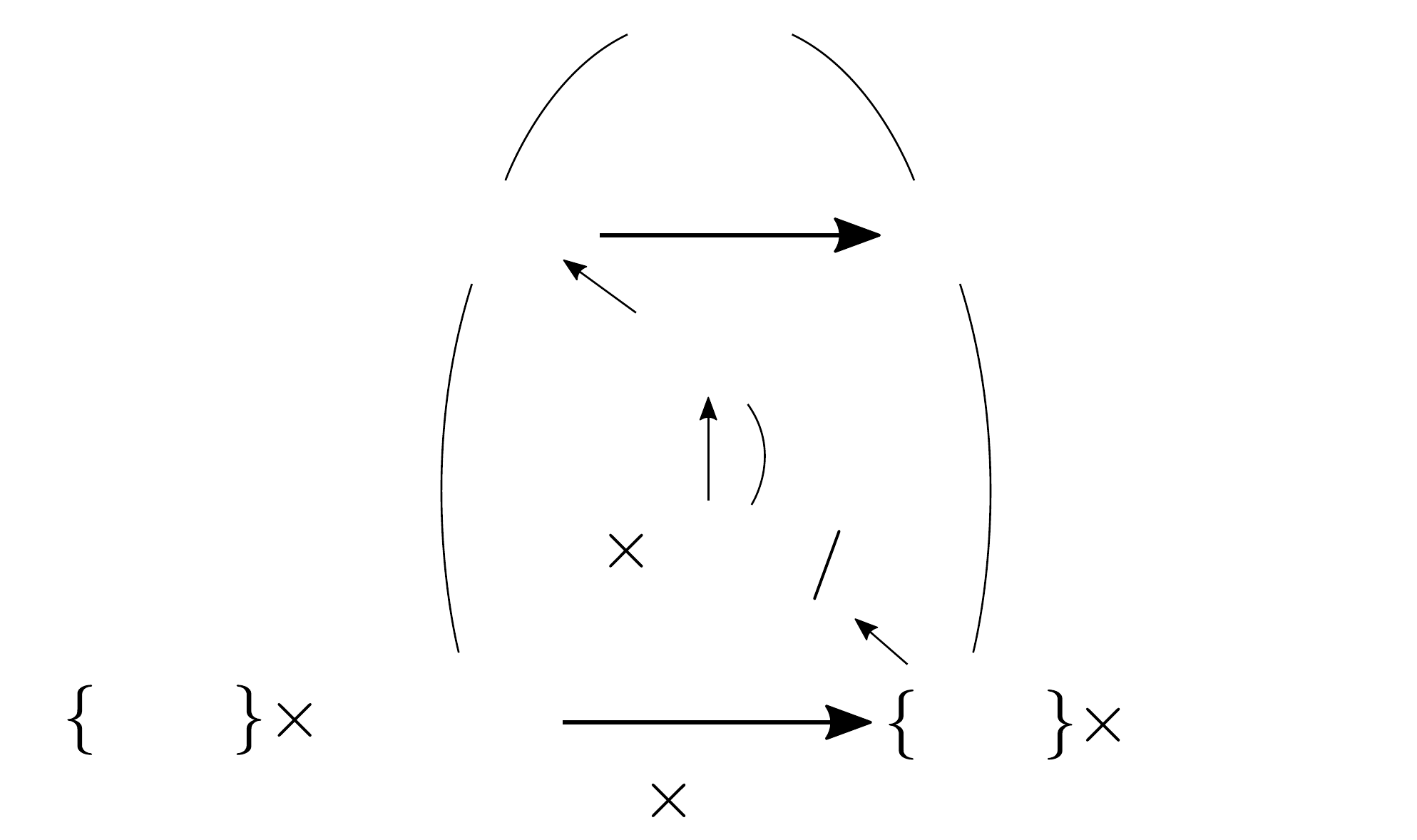
	\endgroup
	\caption{Objects in a computer program involved in an electromagnetic transformation. The blue bits tag {\color{blue}number data used when requesting particular physical field numbers} e.g. for computing a numerical quadrature. The green bits tag predominant {\color{green}computations and fields} necessary to be represented in the computers memory. It is important to note, that an electromagnetic field quantity is represented as one of the green bits, even though its contained degrees of freedom are thought of being blue bits at first. }
	\label{fig:interpreteddiagram}
\end{figure*}

In Sec. \ref{sec:theory} a family of datatypes for well-formed terms in an augmented $\lambda$-calculus was set up. 
That family was indexed by a ``Type'' being the term's type and a ``Context'', to realize a valid use of variables.
This should serve as a foundation for an implementation of our application from Sec. \ref{sec:problem}.
Objects and equivalences from that application, the partial derivative and the multivariate chain rule, can be expressed within this $\lambda$-calculus.
But while guaranteeing these translations to be well-formed terms, this still does not make a computation.
In Sec. \ref{sec:softwarecontext} it was motivated how the formulation of the construction of a shape function element within a programming language easily leads to a meta-implementation.
The meta-implementation's purpose is to generate an efficient implementation where the meta-implementation itself should be focused on validity rather than efficiency.
Our approach in Sec. \ref{sec:theory} should enable to achieve a high validity in a meta-implementation.

We are now focusing on how to give the $\lambda$-terms a suitable interpretation:
Our approach has changed the task of giving a direct interpretation to partial derivatives, into the task of giving a direct interpretation to some lower level primitives.
These are: variable access, quantification, substitution, projection, composition and the two \textit{special functions} which were univariate derivative and a scalar multiplication.
Replacing one notion of partial derivative\footnote{or an evaluation of Jacobian matrices if you like so}, with these lot of operations seems quite a lot of machinery and not worth the trade.

At this point we argue that: First, this approach is in some sense minimal.
Decomposing partial derivative into more basic notions as in equation \ref{eq:partial-derivative} involves only a notion of univariate derivative of a certain function\footnote{which is correspondence with what M. Spivak mentioned as \textit{ordinary} derivative and cited in equation \ref{eq:ordinaryderivativeofcertainfunction}}. 
With our elaboration in section \ref{sec:problem} we collected what obligations arise when working out such decomposition in a way, precise enough to reach a level of \textit{``correctness and completeness necessary to get a computer program to work''}\cite{thurston94}. This level might be \textit{``a couple of orders of magnitude higher''}\cite{thurston94} than the level it needs to convince humans.

Second, one might \textit{have} nothing more than coordinate transition functions in an implementation on a machine as representations of the \textit{data} of boundary value problems.
Recall that motivation from section \ref{sec:electromagneticalcontext} which promised one generic rule for coordinate transformations, indexed by three indices $(p,q,\omega)$.
These should cover the electromagnetic quantities of interest that will show up in an implementation, applying techniques mentioned in section \ref{sec:softwarecontext}.
There we have identified as an important ingredient the \textit{degrees of freedom} which were mappings from a polynomial differential form $u$ to a real number:

$$u\mapsto \int_f(\text{tr}_f u)\wedge q : \mathcal{P}_{\bullet}\Lambda^k(T)\to\mathbb{R}.$$

Here we would identify $u$ as the section of an associated bundle as in section \ref{sec:electromagneticalcontext}.
Since this section might not be directly representable within the machine's memory, we would handle a coordinate representation of it instead.
This was denoted $\color{green}{\stackrel{xyz}{\sigma}}$ in Fig. \ref{fig:diagram} or {\color{green}{phys. field numbers corresp. to xyz}} in Fig. \ref{fig:interpreteddiagram}.

Evaluation of these integrals in order to produce data for a discrete linear system to be solved, is usually done by numerical quadrature.
That quadrature \textit{queries} the physical {\color{blue}{field quantitie's value}} at {\color{blue}{some specific coordinates}}.
If one cannot, or does not want to, predict which physical field values will be queried, it might appeal to {\color{green}{represent the physical field as a computation}}.
Of course, that computation might internally interpolate field values at some specific coordinates with polynomials, as is the case with the polynomial differential forms.


When speaking of finite elements, these are usually transformed to a \textit{reference element} already in order to evaluate a numerical quadrature.
Such approach is desired, since it allows the quadrature to be implemented in a fixed way with precalculated coefficents.
Therefore we might say that every computer program implementing this technique has to deal with partial derivatives in some way already.
But having just one integral transformation might not be worth the effort we have made in the previous part of this paper.
We think that a construction of the boundary value problem itself might be given in terms of a longer chain of composed coordinate transformations.
The representation of a boundary value problem could then be internally encoded as an equivalence transformation out of primitives.
But this is just a motivation for our approach.

We argue it to be, at least, a justified perspective that these kinds of coordinate transformations\footnote{or integral transformations if you like so} apply for a wide range of \textit{numerical software} as sketched in section \ref{sec:softwarecontext}.
Here the particular focus was on boundary value problems expressable by some Hodge-Laplacian over a manifold, discretized with a simplicial complex.

In a broader sense, we understand a part of \textit{software} as a mapping of mathematical models for coordinate transformations, given by (\ref{eq:generaltransformationlaw}) in terms of the partial derivative, into a programming language.
This is done by decomposing this mapping into first, a mapping of partial derivative into a formal language based on $\lambda$-calculus, and second, a mapping of these $\lambda$-calculus terms to computations in a concrete programming language.
Benefit arises, since the first mapping can be discussed and justified on a theoretical basis, where the second is much more arbitrary in its nature: arbitrary in a sense that we have to deal with different computational models to create programs that run on different kinds of machines.
Now, we will elaborate on some of these more arbitrary ways to map our specific $\lambda$-terms to a concrete programming language, or rather map them to a model of evaluation coming with such concrete programming language. A concrete programming language for that purpose comes with a syntax \textit{and} an evaluation strategy.


\subsection{Interpreting $\lambda$-terms formally}

A straight-forward way, since we already have an Agda formalization, would be to continue here by implementing an evaluation function.
Doing so is a typical task and two new concepts occur: for a context $\Gamma$ and one of our custom types $a$, the evaluation function ``eval'' maps an \textit{environment} of that context and a $\lambda$-term of that type and context to a \textit{value} of the \textit{interpretation} of $a$.
$$\text{eval}:\text{Environment }\Gamma\to\text{ Term }\Gamma\;a\to\text{ Interpretation }a$$
Here, ``Interpretation'' maps our custom types from Fig. \ref{fig:customtypes} to types of the Agda language which are elements of the universe ``Set'':
$$\text{Interpretation } : \text{ Type }\to\text{ Set}.$$
An interpretation of our $\lambda$-term's types then really is captured by a function from our previously defined ``Type'' to the types of the programming language which is Agda in this case.

We showed in Sec. \ref{sec:theory} that a context - in the way we introduced it - can be regarded as a list holding multiple variable name and type combinations.
An \textit{environment} for such context can be regarded as holding the corresponding values of these types.
We might only use an environment by \textit{looking up variables}, which are de Bruijn indices in our formalization:
$$\text{lookup }:\text{ Environment } \Gamma \to a \in \Gamma \to \text{ Interpretation } a.$$
This shows that the choice of implementing an environment is already a little less \textit{fixed}.
We could mimick the context and use a list, but in contrast to the context, the environment will be present in our evaluation's \textit{computation} where we might forget about the context completely.
One might prove within Agda that an implementation of lookup never fails when given a valid de Bruijn index and then strip all the type information, revealing bare computations.
Therefore with the environment, we do want to incorporate some aspects of performance.
Speaking of performance, we might have a hard time continuing to use the Agda language itself as a target for evaluation.
It is possible to implement a numerical software within this language\footnote{The Agda language is implemented in Haskell and can use Haskell methods, which in turn via a \textit{foreign function interface} can call arbitrary system libraries}, but this programming language's environment offers only a limited help to put the machine into its most efficient state for the purpose of a computation.
A typical goal is to map scalars to \textit{unboxed}\footnote{Programming languages with automatic reference and memory management tend to implicitly attach typing information to a value to be able treating this value via references instead. This is done because references into memory on a machine have a uniform representation. It is called ``boxing'' of a value. Boxing often demands memory allocation. Preventing frequent memory allocation, e.g. per-scalar memory allocation, is very important to achieve high efficiency in an implementation.} floating point machine numbers which lack a lot of the properties of their counterparts from $\mathbb{R}$.
The Agda programming language offers more support for showing properties on rationals or constructive real numbers.
Unfortunately, these numbers tend to have a representation making them unsuited for numerical computations.
But that is a usual trade we have a lot with numerical algorithms: once their theory is worked out for exact real numbers and the algorithm is \textit{stable}, then we do apply it on inexact floating point machine numbers, fingers crossed\footnote{Unless using promising techniques such as interval arithmetic to provide guarantees for this approach}.

What we also get with the evaluation as a mapping is the possibility to proof preservation of the $=_\text{lin}$ and $=_\text{chain}$ equivalences from Sec. \ref{sec:problem} for our intended implementation.
If we chose to interpret our custom $\lambda$-calculus functions really as functions, then in particular the univariate derivative might be chosen to be an inexact black-box operation such as the difference quotient.
In that case, the chances are high that we will loose the possibility to exactly proof that $=_\text{lin}$ and $=_\text{chain}$ are preserved by our evaluation function even if operating on exact rational numbers.
That could be intended and we might formally track error bounds with all our operations to proof that the error resulting from this operation amortizes comparing to some other error.
But the more interesting case would be, to interpret $\lambda$-functions not as functions but as data structures with a more \textit{interesting} interpretation of $\lambda$-function-application as a data transformation.
This perspective is elaborated in section \ref{sec:datatransformations}.

\subsection{Interpreting $\lambda$-terms less formally}

One might have a formal model of a programming language at hand such that $\lambda$-terms can be translated.
Depending on the degree of formalism, surjectivity of the eval function can be proven or even that the resulting terms are still well-defined in the target language.
This is essentially some type of code generation, where the weakest variant would be to interpret all our custom $\lambda$-terms within the string monoid, calling it \textit{``code''}.
Even if one does not formally model the target language, this still gives a possibility to implement transformations at the $\lambda$-term level, before they are evaluated to code.
Although introducing a large margin for interpretation and bugs, this approach could be well suited for generating code running in a very limited, e.g. lock-step, environment.

\subsection{Interpreting $\lambda$-terms as data transformations}\label{sec:datatransformations}

As motivated before, an interpretation of $\lambda$-functions to data structures of a target language might currently be the most rewarding one.
For the finite element spaces from our application, a lot of different multivariate polynomial functions have to be operated with.
These functions can be represented by polynomial coefficients, together with a custom function application operation that does use these coefficients to compute the polynomial.
Furthermore, the univariate derivative operation on such coefficient representation is not only a very cheap one, but also exact when using exact number representations.
That enables to proof the eval function to preserve $=_\text{lin}$ and $=_\text{chain}$ computationally.

In a meta-implementation for generating an efficient implementation it seems resonable to start out with an initial candidate for the implementation and then apply rewrite rules to optimize this implementation.
Implementing correct rewrites and data transformations is where Agda, and functional programming in general, shines.
The reason for that is an inductive definition of the data structures in question which enables an \textit{exhaustion check} to proof functions to be total, i.e. not having missed a case.
This usually pays off in case-analysis-heavy applications such as designing domain specific languages, as we do here, or improving the encoding of a data structure.
As for multivariate polynomials, these can obviously be represented as \textit{packed} chunks of computer memory, holding their coefficients.
But we might add some information or invent an interesting reference type for better tying them to the simplicial complex they are originating from.
Having equivalence proven for one \textit{obvious} encoding it can be easier for a new encoding to show it isomorphic, transferring the proofs.

\section{Conclusion}

We have explained transformations on the partial derivative in terms of computational notions from $\lambda$-calculus
with an additional term substitution. This mechanism has been implemented to generate listings for the general case as in appendix A and for all concrete multivariate cases, indexed by ${\color{blue}\mathbf{j}}~\in~\mathbb{N}$ and ${\color{blue}\mathbf{i}}\in[1,{\color{blue}\mathbf{j}}]$, exemplary for ${\color{blue}\mathbf{j}}=2$ and ${\color{blue}\mathbf{i}}=1$ as in appendix B, out of the same internal representation. 
It was argued, what general obligations arise when translating the theory into a computational layer of abstraction, for which the $\lambda$-calculus served as a model.
We showed how a translation into an augmented $\lambda$-calculus can be formalized within type theoretical terms and implemented that formalization in the Agda programming language.
Finally, we gave some examples how to make use of the presented approach and favorized one particular possibility.
Our current research is about this exact undertaking and the foundational considerations are shown in our contribution.

Small programs as well as big software, no matter whether directly implementing this layer or not, will suffer from the inevitable tediousness of coordinate transformations when exploiting these techniques too much. That does not pose a problem when being aware of this issue and actively increasing rigor if this kind of complexity gets out of control. We have presented a way to establish that direction of rigor, motivated by the application of encoding the transformation laws common to the electromagnetic theory. Accompanying that way is an interpretation to guide an implementation demanding it.

\section{Appendix}\label{sec7}

In the appendix we give a listing of the computational equivalences used to demonstrate the dependencies of the notion of partial derivative and the chain rule of the partial derivative on the notion of univariate derivative and the corresponding univariate chain rule. Both listings have been created out of the same internal representation with the rules of parameter-pack expansion borrowed from the C++ programming language, with the help of our own implementation of the parameter-pack expansion, supporting the mentioned substitution. For the expanded listing in \ref{sec7b} we chose $f, g : \mathbb{R}^2\rightarrow\mathbb{R}^2$ and ${\color{blue}\mathbf{i}}=1$.

Furthermore we attached a translation of the Agda datatype of well-formed $\lambda$-terms from Fig. \ref{fig:agdaterm}.

%
%

\onecolumn

\subsection{Appendix A}\label{sec7a}
\setcounter{equation}{0}
\def\theequation{A\arabic{equation}}
\begingroup
\fontsize{8pt}{10pt}\selectfont
\addtolength{\jot}{-0.4em} 
\begin{align}{}
& &&\frac{\partial \left(f^{\color{blue}\mathbf{j}}\circ g\right)}{\partial {x}^{{\color{blue}\mathbf{i}}}}\label{eq:A1}\\
&\equiv_{\text{def}}&&\left(\underline{x}...\right)\mapsto {\left({\color{purple}\mathbf{z}}\mapsto \left(f^{\color{blue}\mathbf{j}}\circ g\right)\left(\underline{x}...\left[\bullet^{{\color{blue}\mathbf{i}}}:={\color{purple}\mathbf{z}}\right]\right)\right)}^{\color{red}\mathbf{'}}\left({x}^{{\color{blue}\mathbf{i}}}\right)\label{eq:A2}\\
&\equiv_{\circ}&&\left(\underline{x}...\right)\mapsto {\left({\color{purple}\mathbf{z}}\mapsto f^{\color{blue}\mathbf{j}}\left(\underline{g}\left(\underline{x}...\left[\bullet^{{\color{blue}\mathbf{i}}}:={\color{purple}\mathbf{z}}\right]\right)...\right)\right)}^{\color{red}\mathbf{'}}\left({x}^{{\color{blue}\mathbf{i}}}\right)\label{eq:A3}\\
&=_{\text{lin}'}&&\left(\underline{x}...\right)\mapsto {\color{green}\sum_k}{\left({\color{purple}\mathbf{z}}\mapsto f^{\color{blue}\mathbf{j}}\left(\underline{g}\left(\underline{x}...\right)...\left[\bullet^{{\color{green}k}} :=g^{\color{green}k}\left(\underline{x}...\left[\bullet^{{\color{blue}\mathbf{i}}}:={\color{purple}\mathbf{z}}\right]\right)\right]\right)\right)}^{\color{red}\mathbf{'}}\left({x}^{{\color{blue}\mathbf{i}}}\right)\label{eq:A4}\\
&\equiv_{\circ}&&\left(\underline{x}...\right)\mapsto {\color{green}\sum_k}{\left(\overbrace{\left({\color{purple}\mathbf{z}}\mapsto f^{\color{blue}\mathbf{j}}\left(\underline{g}\left(\underline{x}...\right)...\left[\bullet^{{\color{green}k}} :={\color{purple}\mathbf{z}}\right]\right)\right)}^{{\color{green}a^{\color{green}k}}}\circ \overbrace{\left({\color{purple}\mathbf{z}}\mapsto g^{\color{green}k}\left(\underline{x}...\left[\bullet^{{\color{blue}\mathbf{i}}}:={\color{purple}\mathbf{z}}\right]\right)\right)}^{{\color{green}b^{\color{green}k}}}\right)}^{\color{red}\mathbf{'}}\left({x}^{{\color{blue}\mathbf{i}}}\right)\label{eq:A5}\\
&\equiv_{\text{def}}&&\left(\underline{x}...\right)\mapsto {\color{green}\sum_k}{\left({\color{green}a^{\color{green}k}}\circ {\color{green}b^{\color{green}k}}\right)}^{\color{red}\mathbf{'}}\left({x}^{{\color{blue}\mathbf{i}}}\right)\label{eq:A6}\\
&=_{\text{chain}'}&&\left(\underline{x}...\right)\mapsto {\color{green}\sum_k}\left(\left({{\color{green}a^{\color{green}k}}}^{\color{red}\mathbf{'}}\circ {\color{green}b^{\color{green}k}}\right)\left({x}^{{\color{blue}\mathbf{i}}}\right)\right)\cdot {{\color{green}b^{\color{green}k}}}^{\color{red}\mathbf{'}}\left({x}^{{\color{blue}\mathbf{i}}}\right)\label{eq:A7}\\
&\equiv_{\text{def}}&&\left(\underline{x}...\right)\mapsto {\color{green}\sum_k}\underbrace{\left(\left(\overbrace{{\left({\color{purple}\mathbf{z}}\mapsto f^{\color{blue}\mathbf{j}}\left(\underline{g}\left(\underline{x}...\right)...\left[\bullet^{{\color{green}k}} :={\color{purple}\mathbf{z}}\right]\right)\right)}^{\color{red}\mathbf{'}}}^{{{\color{green}a^{\color{green}k}}}^{\color{red}\mathbf{'}}}\circ \overbrace{\left({\color{purple}\mathbf{z}}\mapsto g^{\color{green}k}\left(\underline{x}...\left[\bullet^{{\color{blue}\mathbf{i}}}:={\color{purple}\mathbf{z}}\right]\right)\right)}^{{\color{green}b^{\color{green}k}}}\right)\left({x}^{{\color{blue}\mathbf{i}}}\right)\right)}_{\frac{\partial f^{\color{blue}\mathbf{j}}}{\partial y^{\color{green}k}}\left(g^{\color{green}k}\left(\underline{x}...\right)\right)}\cdot \underbrace{\overbrace{{\left({\color{purple}\mathbf{z}}\mapsto g^{\color{green}k}\left(\underline{x}...\left[\bullet^{{\color{blue}\mathbf{i}}}:={\color{purple}\mathbf{z}}\right]\right)\right)}^{\color{red}\mathbf{'}}\left({x}^{{\color{blue}\mathbf{i}}}\right)}^{{{\color{green}b^{\color{green}k}}}^{\color{red}\mathbf{'}}\left({x}^{{\color{blue}\mathbf{i}}}\right)}}_{\frac{\partial g^{\color{green}k}}{\partial {x}^{{\color{blue}\mathbf{i}}}}\left(\underline{x}...\right)}\label{eq:A8}\\
&\equiv_{\beta}&&\left(\underline{x}...\right)\mapsto {\color{green}\sum_k}\left(\left({\left({\color{purple}\mathbf{z}}\mapsto f^{\color{blue}\mathbf{j}}\left(\underline{g}\left(\underline{x}...\right)...\left[\bullet^{{\color{green}k}} :={\color{purple}\mathbf{z}}\right]\right)\right)}^{\color{red}\mathbf{'}}\circ \left({\color{purple}\mathbf{z}}\mapsto g^{\color{green}k}\left(\underline{x}...\left[\bullet^{{\color{blue}\mathbf{i}}}:={\color{purple}\mathbf{z}}\right]\right)\right)\right)\left({x}^{{\color{blue}\mathbf{i}}}\right)\right)\cdot \underbrace{\left(\left(\underline{y}...\right)\mapsto {\left({\color{purple}\mathbf{z}}\mapsto g^{\color{green}k}\left(\underline{y}...\left[\bullet^{{\color{blue}\mathbf{i}}}:={\color{purple}\mathbf{z}}\right]\right)\right)}^{\color{red}\mathbf{'}}\left({y}^{{\color{blue}\mathbf{i}}}\right)\right)}_{\frac{\partial g^{\color{green}k}}{\partial {x}^{{\color{blue}\mathbf{i}}}}}\left(\underline{x}...\right)\label{eq:A9}\\
&\equiv_{\text{def}}&&\left(\underline{x}...\right)\mapsto {\color{green}\sum_k}\underbrace{\left(\left({\left({\color{purple}\mathbf{z}}\mapsto f^{\color{blue}\mathbf{j}}\left(\underline{g}\left(\underline{x}...\right)...\left[\bullet^{{\color{green}k}} :={\color{purple}\mathbf{z}}\right]\right)\right)}^{\color{red}\mathbf{'}}\circ \left({\color{purple}\mathbf{z}}\mapsto g^{\color{green}k}\left(\underline{x}...\left[\bullet^{{\color{blue}\mathbf{i}}}:={\color{purple}\mathbf{z}}\right]\right)\right)\right)\left({x}^{{\color{blue}\mathbf{i}}}\right)\right)}_{\frac{\partial f^{\color{blue}\mathbf{j}}}{\partial y^{\color{green}k}}\left(g^{\color{green}k}\left(\underline{x}...\right)\right)}\cdot \frac{\partial g^{\color{green}k}}{\partial {x}^{{\color{blue}\mathbf{i}}}}\left(\underline{x}...\right)\label{eq:A10}\\
&\equiv_{\circ}&&\left(\underline{x}...\right)\mapsto {\color{green}\sum_k}\underbrace{\left({\left({\color{purple}\mathbf{z}}\mapsto f^{\color{blue}\mathbf{j}}\left(\underline{g}\left(\underline{x}...\right)...\left[\bullet^{{\color{green}k}} :={\color{purple}\mathbf{z}}\right]\right)\right)}^{\color{red}\mathbf{'}}\left(g^{\color{green}k}\left(\underline{x}...\right)\right)\right)}_{\frac{\partial f^{\color{blue}\mathbf{j}}}{\partial y^{\color{green}k}}\left(g^{\color{green}k}\left(\underline{x}...\right)\right)}\cdot \frac{\partial g^{\color{green}k}}{\partial {x}^{{\color{blue}\mathbf{i}}}}\left(\underline{x}...\right)\label{eq:A11}\\
&\equiv_{\beta}&&\left(\underline{x}...\right)\mapsto {\color{green}\sum_k}\underbrace{\left(\left(\underline{y}...\right)\mapsto {\left({\color{purple}\mathbf{z}}\mapsto f^{\color{blue}\mathbf{j}}\left(\underline{y}...\left[\bullet^{{\color{green}k}} :={\color{purple}\mathbf{z}}\right]\right)\right)}^{\color{red}\mathbf{'}}\left(y^{\color{green}k}\right)\right)}_{\frac{\partial f^{\color{blue}\mathbf{j}}}{\partial y^{\color{green}k}}}\left(\underline{g}\left(\underline{x}...\right)...\right)\cdot \frac{\partial g^{\color{green}k}}{\partial {x}^{{\color{blue}\mathbf{i}}}}\left(\underline{x}...\right)\label{eq:A12}\\
&\equiv_{\text{def}}&&\left(\underline{x}...\right)\mapsto {\color{green}\sum_k}\frac{\partial f^{\color{blue}\mathbf{j}}}{\partial y^{\color{green}k}}\left(\underline{g}\left(\underline{x}...\right)...\right)\cdot \frac{\partial g^{\color{green}k}}{\partial {x}^{{\color{blue}\mathbf{i}}}}\left(\underline{x}...\right)\label{eq:A13}\\
&\equiv_{\circ}&&\left(\underline{x}...\right)\mapsto {\color{green}\sum_k}\left(\frac{\partial f^{\color{blue}\mathbf{j}}}{\partial y^{\color{green}k}}\circ g\right)\left(\underline{x}...\right)\cdot \frac{\partial g^{\color{green}k}}{\partial {x}^{{\color{blue}\mathbf{i}}}}\left(\underline{x}...\right)\label{eq:A14}\\
&\equiv_{\text{def}}&&\left(\underline{x}...\right)\mapsto {\color{green}\sum_k}\left(\left(\frac{\partial f^{\color{blue}\mathbf{j}}}{\partial y^{\color{green}k}}\circ g\right)\otimes\frac{\partial g^{\color{green}k}}{\partial {x}^{{\color{blue}\mathbf{i}}}}\right)\left(\underline{x}...\right)\label{eq:A15}\\
&\equiv_{\text{def}}&&\left(\underline{x}...\right)\mapsto {\color{green}\sum_k}\left(\frac{\partial f^{\color{blue}\mathbf{j}}}{\partial y^{\color{green}k}}\quad\otimes^{(g\,\times\,\text{id})}\;\frac{\partial g^{\color{green}k}}{\partial {x}^{{\color{blue}\mathbf{i}}}}\right)\left(\underline{x}...\right)\label{eq:A16}\\
&\equiv_{\text{def}}&&{\color{green}\bigoplus_k}\left(\frac{\partial f^{\color{blue}\mathbf{j}}}{\partial y^{\color{green}k}}\quad\otimes^{(g\,\times\,\text{id})}\;\frac{\partial g^{\color{green}k}}{\partial {x}^{{\color{blue}\mathbf{i}}}}\right)\label{eq:A17}\\
&\cong_{\text{T}}&&\phantom{\color{green}\sum_k}\left(\text{J}^{{{\color{blue}\mathbf{j}}}}_{{{\color{green}k}}'}\quad\quad\quad\quad\quad\quad\text{J}^{{{\color{green}k}}'}_{{{\color{blue}\mathbf{i}}}''}\right)\label{eq:A18}\\
&\equiv&&\phantom{\color{green}\sum_k}\text{J}^{{{\color{blue}\mathbf{j}}}}_{{{\color{green}k}}'}\text{J}^{{{\color{green}k}}'}_{{{\color{blue}\mathbf{i}}}''}\label{eq:A19}
\end{align}
\endgroup
\clearpage
\subsection{Appendix B}\label{sec7b}

\setcounter{equation}{0}
\def\theequation{B\arabic{equation}}

\begingroup
\fontsize{8pt}{10pt}\selectfont
%
\begin{align}{}
& &&\frac{\partial \left(f^{\color{blue}\mathbf{j}}\circ g\right)}{\partial {x}^{1}}\label{eq:B1}\\
&\equiv_{\text{def}}&&\left(x^1,x^2\right)\mapsto {\left({\color{purple}\mathbf{z}}\mapsto \left(f^{\color{blue}\mathbf{j}}\circ g\right)\left({\color{purple}\mathbf{z}},x^2\right)\right)}^{\color{red}\mathbf{'}}\left({x}^{1}\right)\label{eq:B2}\\
&\equiv_{\circ}&&\left(x^1,x^2\right)\mapsto {\left({\color{purple}\mathbf{z}}\mapsto f^{\color{blue}\mathbf{j}}\left(g^1\left({\color{purple}\mathbf{z}},x^2\right),g^2\left({\color{purple}\mathbf{z}},x^2\right)\right)\right)}^{\color{red}\mathbf{'}}\left({x}^{1}\right)\label{eq:B3}\\
&=_{\text{lin}'}&&\left(x^1,x^2\right)\mapsto {\left({\color{purple}\mathbf{z}}\mapsto f^{\color{blue}\mathbf{j}}\left(g^1\left({\color{purple}\mathbf{z}},x^2\right),g^2\left(x^1,x^2\right)\right)\right)}^{\color{red}\mathbf{'}}\left({x}^{1}\right)+{\left({\color{purple}\mathbf{z}}\mapsto f^{\color{blue}\mathbf{j}}\left(g^1\left(x^1,x^2\right),g^2\left({\color{purple}\mathbf{z}},x^2\right)\right)\right)}^{\color{red}\mathbf{'}}\left({x}^{1}\right)\label{eq:B4}\\
&\equiv_{\circ}&&\left(x^1,x^2\right)\mapsto {\left(\overbrace{\left({\color{purple}\mathbf{z}}\mapsto f^{\color{blue}\mathbf{j}}\left({\color{purple}\mathbf{z}},g^2\left(x^1,x^2\right)\right)\right)}^{{\color{green}a^1}}\circ \overbrace{\left({\color{purple}\mathbf{z}}\mapsto g^1\left({\color{purple}\mathbf{z}},x^2\right)\right)}^{{\color{green}b^1}}\right)}^{\color{red}\mathbf{'}}\left({x}^{1}\right)+{\left(\overbrace{\left({\color{purple}\mathbf{z}}\mapsto f^{\color{blue}\mathbf{j}}\left(g^1\left(x^1,x^2\right),{\color{purple}\mathbf{z}}\right)\right)}^{{\color{green}a^2}}\circ \overbrace{\left({\color{purple}\mathbf{z}}\mapsto g^2\left({\color{purple}\mathbf{z}},x^2\right)\right)}^{{\color{green}b^2}}\right)}^{\color{red}\mathbf{'}}\left({x}^{1}\right)\label{eq:B5}\\
&\equiv_{\text{def}}&&\left(x^1,x^2\right)\mapsto {\left({\color{green}a^1}\circ {\color{green}b^1}\right)}^{\color{red}\mathbf{'}}\left({x}^{1}\right)+{\left({\color{green}a^2}\circ {\color{green}b^2}\right)}^{\color{red}\mathbf{'}}\left({x}^{1}\right)\label{eq:B6}\\
&=_{\text{chain}'}&&\left(x^1,x^2\right)\mapsto \left(\left({{\color{green}a^1}}^{\color{red}\mathbf{'}}\circ {\color{green}b^1}\right)\left({x}^{1}\right)\right)\cdot {{\color{green}b^1}}^{\color{red}\mathbf{'}}\left({x}^{1}\right)+\left(\left({{\color{green}a^2}}^{\color{red}\mathbf{'}}\circ {\color{green}b^2}\right)\left({x}^{1}\right)\right)\cdot {{\color{green}b^2}}^{\color{red}\mathbf{'}}\left({x}^{1}\right)\label{eq:B7}\\
&\equiv_{\text{def}}&&\left(x^1,x^2\right)\mapsto \underbrace{\left(\left(\overbrace{{\left({\color{purple}\mathbf{z}}\mapsto f^{\color{blue}\mathbf{j}}\left({\color{purple}\mathbf{z}},g^2\left(x^1,x^2\right)\right)\right)}^{\color{red}\mathbf{'}}}^{{{\color{green}a^1}}^{\color{red}\mathbf{'}}}\circ \overbrace{\left({\color{purple}\mathbf{z}}\mapsto g^1\left({\color{purple}\mathbf{z}},x^2\right)\right)}^{{\color{green}b^1}}\right)\left({x}^{1}\right)\right)}_{\frac{\partial f^{\color{blue}\mathbf{j}}}{\partial y^1}\left(g^1\left(x^1,x^2\right)\right)}\cdot \underbrace{\overbrace{{\left({\color{purple}\mathbf{z}}\mapsto g^1\left({\color{purple}\mathbf{z}},x^2\right)\right)}^{\color{red}\mathbf{'}}\left({x}^{1}\right)}^{{{\color{green}b^1}}^{\color{red}\mathbf{'}}\left({x}^{1}\right)}}_{\frac{\partial g^1}{\partial {x}^{1}}\left(x^1,x^2\right)}+\;\text{one more term}\label{eq:B8}\\
&\equiv_{\beta}&&\left(x^1,x^2\right)\mapsto \left(\left({\left({\color{purple}\mathbf{z}}\mapsto f^{\color{blue}\mathbf{j}}\left({\color{purple}\mathbf{z}},g^2\left(x^1,x^2\right)\right)\right)}^{\color{red}\mathbf{'}}\circ \left({\color{purple}\mathbf{z}}\mapsto g^1\left({\color{purple}\mathbf{z}},x^2\right)\right)\right)\left({x}^{1}\right)\right)\cdot \underbrace{\left(\left(y^1,y^2\right)\mapsto {\left({\color{purple}\mathbf{z}}\mapsto g^1\left({\color{purple}\mathbf{z}},y^2\right)\right)}^{\color{red}\mathbf{'}}\left({y}^{1}\right)\right)}_{\frac{\partial g^1}{\partial {x}^{1}}}\left(x^1,x^2\right)+\;\text{one more term}\label{eq:B9}\\
&\equiv_{\text{def}}&&\left(x^1,x^2\right)\mapsto \underbrace{\left(\left({\left({\color{purple}\mathbf{z}}\mapsto f^{\color{blue}\mathbf{j}}\left({\color{purple}\mathbf{z}},g^2\left(x^1,x^2\right)\right)\right)}^{\color{red}\mathbf{'}}\circ \left({\color{purple}\mathbf{z}}\mapsto g^1\left({\color{purple}\mathbf{z}},x^2\right)\right)\right)\left({x}^{1}\right)\right)}_{\frac{\partial f^{\color{blue}\mathbf{j}}}{\partial y^1}\left(g^1\left(x^1,x^2\right)\right)}\cdot \frac{\partial g^1}{\partial {x}^{1}}\left(x^1,x^2\right)+\;\text{one more term}\label{eq:B10}\\
&\equiv_{\circ}&&\left(x^1,x^2\right)\mapsto \underbrace{\left({\left({\color{purple}\mathbf{z}}\mapsto f^{\color{blue}\mathbf{j}}\left({\color{purple}\mathbf{z}},g^2\left(x^1,x^2\right)\right)\right)}^{\color{red}\mathbf{'}}\left(g^1\left(x^1,x^2\right)\right)\right)}_{\frac{\partial f^{\color{blue}\mathbf{j}}}{\partial y^1}\left(g^1\left(x^1,x^2\right)\right)}\cdot \frac{\partial g^1}{\partial {x}^{1}}\left(x^1,x^2\right)+\;\text{one more term}\label{eq:B11}\\
&\equiv_{\beta}&&\left(x^1,x^2\right)\mapsto \underbrace{\left(\left(y^1,y^2\right)\mapsto {\left({\color{purple}\mathbf{z}}\mapsto f^{\color{blue}\mathbf{j}}\left({\color{purple}\mathbf{z}},y^2\right)\right)}^{\color{red}\mathbf{'}}\left(y^1\right)\right)}_{\frac{\partial f^{\color{blue}\mathbf{j}}}{\partial y^1}}\left(g^1\left(x^1,x^2\right),g^2\left(x^1,x^2\right)\right)\cdot \frac{\partial g^1}{\partial {x}^{1}}\left(x^1,x^2\right)+\;\text{one more term}\label{eq:B12}\\
&\equiv_{\text{def}}&&\left(x^1,x^2\right)\mapsto \frac{\partial f^{\color{blue}\mathbf{j}}}{\partial y^1}\left(g^1\left(x^1,x^2\right),g^2\left(x^1,x^2\right)\right)\cdot \frac{\partial g^1}{\partial {x}^{1}}\left(x^1,x^2\right)+\frac{\partial f^{\color{blue}\mathbf{j}}}{\partial y^2}\left(g^1\left(x^1,x^2\right),g^2\left(x^1,x^2\right)\right)\cdot \frac{\partial g^2}{\partial {x}^{1}}\left(x^1,x^2\right)\label{eq:B13}\\
&\equiv_{\circ}&&\left(x^1,x^2\right)\mapsto \left(\frac{\partial f^{\color{blue}\mathbf{j}}}{\partial y^1}\circ g\right)\left(x^1,x^2\right)\cdot \frac{\partial g^1}{\partial {x}^{1}}\left(x^1,x^2\right)+\left(\frac{\partial f^{\color{blue}\mathbf{j}}}{\partial y^2}\circ g\right)\left(x^1,x^2\right)\cdot \frac{\partial g^2}{\partial {x}^{1}}\left(x^1,x^2\right)\label{eq:B14}\\
&\equiv_{\text{def}}&&\left(x^1,x^2\right)\mapsto \left(\left(\frac{\partial f^{\color{blue}\mathbf{j}}}{\partial y^1}\circ g\right)\otimes\frac{\partial g^1}{\partial {x}^{1}}\right)\left(x^1,x^2\right)+\left(\left(\frac{\partial f^{\color{blue}\mathbf{j}}}{\partial y^2}\circ g\right)\otimes\frac{\partial g^2}{\partial {x}^{1}}\right)\left(x^1,x^2\right)\label{eq:B15}\\
&\equiv_{\text{def}}&&\left(x^1,x^2\right)\mapsto \left(\frac{\partial f^{\color{blue}\mathbf{j}}}{\partial y^1}\quad\otimes^{(g\,\times\,\text{id})}\;\frac{\partial g^1}{\partial {x}^{1}}\right)\left(x^1,x^2\right)+\left(\frac{\partial f^{\color{blue}\mathbf{j}}}{\partial y^2}\quad\otimes^{(g\,\times\,\text{id})}\;\frac{\partial g^2}{\partial {x}^{1}}\right)\left(x^1,x^2\right)\label{eq:B16}\\
&\equiv_{\text{def}}&&\left(\frac{\partial f^{\color{blue}\mathbf{j}}}{\partial y^1}\quad\otimes^{(g\,\times\,\text{id})}\;\frac{\partial g^1}{\partial {x}^{1}}\right)\oplus\left(\frac{\partial f^{\color{blue}\mathbf{j}}}{\partial y^2}\quad\otimes^{(g\,\times\,\text{id})}\;\frac{\partial g^2}{\partial {x}^{1}}\right)\label{eq:B17}\\
&\cong_{\text{T}}&&\left(\text{J}^{{{\color{blue}\mathbf{j}}}}_{{1}'}\quad\quad\quad\quad\quad\quad\text{J}^{{1}'}_{{1}''}\right)+_{T}\left(\text{J}^{{{\color{blue}\mathbf{j}}}}_{{2}'}\quad\quad\quad\quad\quad\quad\text{J}^{{2}'}_{{1}''}\right)\label{eq:B18}\\
&\equiv&&\text{J}^{{{\color{blue}\mathbf{j}}}}_{{1}'}\text{J}^{{1}'}_{{1}''}+_{T}\text{J}^{{{\color{blue}\mathbf{j}}}}_{{2}'}\text{J}^{{2}'}_{{1}''}\label{eq:B19}
\end{align}
\endgroup
\clearpage

\twocolumn

\subsection{Appendix C}\label{sec:allrules}

\setcounter{equation}{0}
\def\theequation{C\arabic{equation}}

\begin{prooftree}
\AxiomC{$i : (x , a) \in \Gamma$}
\RightLabel{($- \because - $)}
\inference[1]{$\Gamma \vdash x \;\because\; i : a$}
\end{prooftree}

\begin{prooftree}
\AxiomC{$\Gamma , (x , \tuple m) \vdash T : \tuple n$}
\RightLabel{($\lambda^t - {\;}^t - {\;.\;} -$)}
\inference[1]{$\Gamma \vdash \lambda^t m {\;}^t x {\;.\;} T : \funMN m\;n$}
\end{prooftree}

\begin{prooftree}
\AxiomC{$\Gamma , (x , \tuple m) \vdash T : \scalar $}
\RightLabel{($\lambda^t - {\;}^{s} - {\;.\;} -$)}
\inference[1]{$\Gamma \vdash \lambda^t m {\;}^{s} x {\;.\;} T : \funMI m$}
\end{prooftree}

\begin{prooftree}
\AxiomC{$\Gamma , (x , \scalar ) \vdash T : \tuple n$}
\RightLabel{($\lambda^{st} - {\;.\;} -$)}
\inference[1]{$\Gamma \vdash \lambda^{st} x {\;.\;} T : \funIN n$}
\end{prooftree}

\begin{prooftree}
\AxiomC{$\Gamma , (x , \scalar ) \vdash T : \scalar $}
\RightLabel{($\lambda^{ss} - {\;.\;} -$)}
\inference[1]{$\Gamma \vdash \lambda^{ss} x {\;.\;} T : \funII$}
\end{prooftree}

\begin{prooftree}
\AxiomC{$\Gamma , (x , \Index k) \vdash T : \scalar $}
\RightLabel{($\Sigma - {\;}^{s} - {\;.\;} -$)}
\inference[1]{$\Gamma \vdash \Sigma\; k {\;}^{s} x {\;.\;} T : \scalar$}
\end{prooftree}

\begin{prooftree}
\AxiomC{$\Gamma \vdash T : \funII$}
\AxiomC{$\Gamma \vdash U : \scalar $}
\RightLabel{($- {\;}^{11} \llparenthesis - \rrparenthesis$)}
\inference[2]{$\Gamma \vdash T {\;}^{11} \llparenthesis\; U \;\rrparenthesis : \scalar$}
\end{prooftree}

\begin{prooftree}
\AxiomC{$\Gamma \vdash T : \funMI m $}
\AxiomC{$\Gamma \vdash U : \tuple m$}
\RightLabel{($- {\;}^{m1} \llparenthesis - \rrparenthesis$)}
\inference[2]{$\Gamma \vdash T {\;}^{m1} \llparenthesis\; U \;\rrparenthesis : \scalar$}
\end{prooftree}

\begin{prooftree}
\AxiomC{$\Gamma \vdash T : \funIN n $}
\AxiomC{$\Gamma \vdash U : \scalar $}
\RightLabel{($- {\;}^{1n} \llparenthesis - \rrparenthesis$)}
\inference[2]{$\Gamma \vdash T {\;}^{1n} \llparenthesis\; U \;\rrparenthesis : \tuple n$}
\end{prooftree}

\begin{prooftree}
\AxiomC{$\Gamma \vdash T : \funMN m\;n$}
\AxiomC{$\Gamma \vdash U : \tuple m$}
\RightLabel{($- {\;}^{mn} \llparenthesis - \rrparenthesis$)}
\inference[2]{$\Gamma \vdash T {\;}^{mn} \llparenthesis\; U \;\rrparenthesis : \tuple n$}
\end{prooftree}

\begin{prooftree}
\AxiomC{$\Gamma \vdash T : \tuple k $}
\AxiomC{$i : \Fin k $}
\AxiomC{$\Gamma \vdash U : \scalar$}
\inference[3]{$\Gamma \vdash T {\;}^{k}[\bullet i := U] : \tuple k$}
\end{prooftree}

\begin{prooftree}
\AxiomC{$\Gamma \vdash T : \funIN n $}
\AxiomC{$i : \Fin n $}
\AxiomC{$\Gamma \vdash U : \scalar$}
\inference[3]{$\Gamma \vdash T {\;}^{n}[\bullet i := U] : \funIN n$}
\end{prooftree}

\begin{prooftree}
\AxiomC{$\Gamma \vdash T : \funMN m\;n$}
\AxiomC{$i : \Fin n $}
\AxiomC{$\Gamma \vdash U : \scalar$}
\inference[3]{$\Gamma \vdash T {\;}^{*}[\bullet i := U] : \funMN m\;n$}
\end{prooftree}

\begin{prooftree}
\AxiomC{$\Gamma \vdash T : \tuple k $}
\AxiomC{$\Gamma \vdash U : \Index k$}
\AxiomC{$\Gamma \vdash V : \scalar$}
\inference[3]{$\Gamma \vdash T {\;}^{ki}[\bullet U := V] : \tuple k$}
\end{prooftree}

\begin{prooftree}
\AxiomC{$\Gamma \vdash T : \funIN n $}
\AxiomC{$\Gamma \vdash U : \Index n$}
\AxiomC{$\Gamma \vdash V : \scalar$}
\inference[3]{$\Gamma \vdash T {\;}^{ni}[\bullet U := V] : \funIN n$}
\end{prooftree}

\begin{prooftree}
\AxiomC{$\Gamma \vdash T : \funMN m\;n$}
\AxiomC{$\Gamma \vdash U : \Index n$}
\AxiomC{$\Gamma \vdash V : \scalar$}
\inference[3]{$\Gamma \vdash T {\;}^{*i}[\bullet U := V] : \funMN m\;n$}
\end{prooftree}

\begin{prooftree}
\AxiomC{$\Gamma \vdash T : \tuple k $}
\AxiomC{$i : \Fin k $}
\RightLabel{($- \;\hat{}{\,}^{k} -$)}
\inference[2]{$\Gamma \vdash T \;\hat{}{\,}^{k} i : \scalar$}
\end{prooftree}

\begin{prooftree}
\AxiomC{$\Gamma \vdash T : \funIN n $}
\AxiomC{$i : \Fin n $}
\RightLabel{($- \;\hat{}{\,}^{n} -$)}
\inference[2]{$\Gamma \vdash T \;\hat{}{\,}^{n} i : \funII$}
\end{prooftree}

\begin{prooftree}
\AxiomC{$\Gamma \vdash T : \funMN m\;n$}
\AxiomC{$i : \Fin n $}
\RightLabel{($- \;\hat{}{\,}^{*} -$)}
\inference[2]{$\Gamma \vdash T \;\hat{}{\,}^{*} i : \funMI m$}
\end{prooftree}

\begin{prooftree}
\AxiomC{$\Gamma \vdash T : \tuple k $}
\AxiomC{$\Gamma \vdash U : \Index k$}
\RightLabel{($- \;\hat{}{\,}^{ki} -$)}
\inference[2]{$\Gamma \vdash T \;\hat{}{\,}^{ki} U : \scalar$}
\end{prooftree}

\begin{prooftree}
\AxiomC{$\Gamma \vdash T : \funIN n $}
\AxiomC{$\Gamma \vdash U : \Index n$}
\RightLabel{($- \;\hat{}{\,}^{ni} -$)}
\inference[2]{$\Gamma \vdash T \;\hat{}{\,}^{ni} U : \funII$}
\end{prooftree}

\begin{prooftree}
\AxiomC{$\Gamma \vdash T : \funMN m\;n$}
\AxiomC{$\Gamma \vdash U : \Index n$}
\RightLabel{($- \;\hat{}{\,}^{*i} -$)}
\inference[2]{$\Gamma \vdash T \;\hat{}{\,}^{*i} U : \funMI m$}
\end{prooftree}

\begin{prooftree}
\AxiomC{$\Gamma \vdash T : \funII $}
\AxiomC{$\Gamma \vdash U : \funII $}
\RightLabel{($- \circ^{111} -$)}
\inference[2]{$\Gamma \vdash T \circ^{111} U : \funII$}
\end{prooftree}

\begin{prooftree}
\AxiomC{$\Gamma \vdash T : \funMI k $}
\AxiomC{$\Gamma \vdash U : \funIN k$}
\RightLabel{($- \circ^{1k1} -$)}
\inference[2]{$\Gamma \vdash T \circ^{1k1} U : \funII$}
\end{prooftree}

\begin{prooftree}
\AxiomC{$\Gamma \vdash T : \funMN k n$}
\AxiomC{$\Gamma \vdash U : \funMN m k$}
\RightLabel{($- \circ^{mkn} -$)}
\inference[2]{$\Gamma \vdash T \circ^{mkn} U : \funMN m\;n$}
\end{prooftree}

\begin{prooftree}
\AxiomC{$\Gamma \vdash T : \funIN n$}
\AxiomC{$\Gamma \vdash U : \funMI m $}
\RightLabel{($- \circ^{m1n} -$)}
\inference[2]{$\Gamma \vdash T \circ^{m1n} U : \funMN m\;n$}
\end{prooftree}

\begin{prooftree}
\AxiomC{$\Gamma \vdash T : \funII$}
\RightLabel{($- {\;}^{'s}$)}
\inference[1]{$\Gamma \vdash T {\;}^{'s} : \funII$}
\end{prooftree}

\begin{prooftree}
\AxiomC{$\Gamma \vdash T : \scalar$}
\AxiomC{$\Gamma \vdash U : \scalar$}
\RightLabel{($- \cdot^{s} -$)}
\inference[2]{$\Gamma \vdash T \cdot^{s} U : \scalar$}
\end{prooftree}

\bibliographystyle{abbrv}
\bibliography{wileyNJD-AMA}

\begin{thebibliography}{10}

\bibitem{fenics}
M.~S. Aln{\ae}s, J.~Blechta, J.~Hake, A.~Johansson, B.~Kehlet, A.~Logg,
  C.~Richardson, J.~Ring, M.~E. Rognes, and G.~N. Wells.
\newblock The fenics project version 1.5.
\newblock {\em Archive of Numerical Software}, 3(100), 2015.

\bibitem{arnold09}
D.~N. Arnold, R.~S. Falk, and R.~Winther.
\newblock Geometric decompositions and local bases for spaces of finite element
  differential forms.
\newblock {\em Comput. Methods Appl. Mech. Engrg.}, 2009.

\bibitem{arnold10}
D.~N. Arnold, R.~S. Falk, and R.~Winther.
\newblock Finite element exterior calculus: from hodge theory to numerical
  stability.
\newblock {\em Bull. Amer. Math. Soc. 47 (2010), 281-354}, 2010.

\bibitem{auchmann14}
B.~Auchmann and S.~Kurz.
\newblock Observers and splitting structures in relativistic electrodynamics.
\newblock {\em Journal of Physics A: Mathematical and Theoretical}, 47:435202,
  10 2014.

\bibitem{baez94}
J.~Baez and J.~P. Muniain.
\newblock {\em Gauge Fields, Knots and Gravity}.
\newblock World Scientific Publishing Company, 1994.

\bibitem{barendregt85}
H.~P. Barendregt.
\newblock {\em The lambda calculus : its syntax and semantics}.
\newblock North-Holland Pub. Co. ; sole distributors for the U.S.A. and Canada
  Elsevier North-Holland Amsterdam ; New York : New York, 1981.

\bibitem{bossavit98}
A.~Bossavit.
\newblock Appendix a - mathematical background.
\newblock In A.~Bossavit, editor, {\em Computational Electromagnetism},
  Electromagnetism, pages 263 -- 318. Academic Press, San Diego, 1998.

\bibitem{Bossavit2001}
A.~Bossavit.
\newblock On the notion of anisotropy of constitutive laws: Some implications
  of the 'hodge implies metric' result.
\newblock {\em COMPEL - The international journal for computation and
  mathematics in electrical and electronic engineering}, 20(1):233--239, 2001.

\bibitem{Bossavit2012}
A.~Bossavit.
\newblock The premetric approach to electromagnetism in the 'waves are not
  vectors' debate.
\newblock {\em Advanced Electromagnetics, 1(1), 97-102}, 2012.

\bibitem{church36}
A.~Church.
\newblock An unsolvable problem of elementary number theory.
\newblock {\em American Journal of Mathematics}, 1936.

\bibitem{debruijn72}
N.~de~Bruijn.
\newblock Lambda calculus notation with nameless dummies, a tool for automatic
  formula manipulation, with application to the church-rosser theorem.
\newblock {\em Indagationes Mathematicae (Proceedings)}, 75(5):381 -- 392,
  1972.

\bibitem{hehl03}
F.~W. Hehl and Y.~N. Obukhov.
\newblock {\em Foundations of Classical Electrodynamics}.
\newblock Birkhauser Boston, 2003.

\bibitem{cpp17}
{ISO}.
\newblock {\em {ISO\slash IEC 14882:2017 Information technology --- Programming
  languages --- C++}}.
\newblock Fifth edition, Dec. 2017.

\bibitem{martinlof83}
P.~Martin-L\"of.
\newblock On the meanings of the logical constant and the justifications of the
  logical laws. technical report 2.
\newblock Technical report, Scuola di Specializziazione in Logica Matematica,
  Universit\`{a} di Siena.

\bibitem{martinlof85}
P.~Martin-L\"{o}f.
\newblock Constructive mathematics and computer programming.
\newblock In {\em Proc. Of a Discussion Meeting of the Royal Society of London
  on Mathematical Logic and Programming Languages}, pages 167--184, Upper
  Saddle River, NJ, USA, 1985. Prentice-Hall, Inc.

\bibitem{norell2009}
U.~Norell.
\newblock Dependently typed programming in agda.
\newblock {\em Koopman P., Plasmeijer R., Swierstra D. (eds) Advanced
  Functional Programming. AFP 2008. Lecture Notes in Computer Science, vol
  5832. Springer, Berlin, Heidelberg}, 2009.

\bibitem{pellikka2013}
M.~Pellikka, T.~Tarhasaari, S.~Suuriniemi, and L.~Kettunen.
\newblock A programming interface to the riemannian manifold in a finite
  element environment.
\newblock {\em Journal of Computational and Applied Mathematics}, 246:225--233,
  07 2013.

\bibitem{Raumonen2009}
P.~Raumonen.
\newblock {\em Mathematical structures for dimensional reduction and
  equivalence classification of electromagnetic boundary value problems}.
\newblock Tampere University of Technology. Publication~. Tampere University of
  Technology, 6 2009.

\bibitem{spivak71}
M.~Spivak.
\newblock {\em Calculus On Manifolds}.
\newblock W. A. Benjamin; New York, 1965.

\bibitem{thurston94}
W.~P. Thurston.
\newblock On proof and progress in mathematics.
\newblock {\em Bull. Amer. Math. Soc. (N.S.) 30 (1994) 161-177}, 1994.

\bibitem{Tonti2013}
E.~Tonti.
\newblock {\em The Mathematical Structure of Classical and Relativistic
  Physics}.
\newblock Birkh\"auser Basel, 2013.

\bibitem{VanDantzig1934}
D.~van Dantzig.
\newblock The fundamental equations of electromagnetism, independent of
  metrical geometry.
\newblock {\em Mathematical Proceedings of the Cambridge Philosophical
  Society}, 30(4):421--427, 1934.

\bibitem{wadler15}
P.~Wadler.
\newblock Propositions as types.
\newblock {\em Commun. ACM}, 58(12):75--84, Nov. 2015.

\bibitem{mathematica}
{Wolfram Research Inc.}
\newblock {\em Mathematica 11.0~}, 2018.

\end{thebibliography}
\end{document}